\documentclass [12pt]{amsart}
\usepackage{times}
\usepackage{amssymb,latexsym,amsmath,amsfonts, eucal}
\usepackage[usenames,dvipsnames,svgnames,table]{xcolor}

\newcommand\hlight[1]{\tikz[overlay, remember picture,baseline=-\the\dimexpr\fontdimen22\textfont2\relax]\node[rectangle,fill=white!50,rounded corners,fill opacity = 0.2,draw,thick,text opacity =1] {$#1$};}

\usepackage[bookmarksnumbered, colorlinks, plainpages]{hyperref}

\newtheorem{thm}{Theorem}[section]
\newtheorem{cor}[thm]{Corollary}
\newtheorem{lem}[thm]{Lemma}

\theoremstyle{definition}
\newtheorem{defn}[thm]{Definition}
\newtheorem{rem}[thm]{Remark}
\numberwithin{equation}{section}

\newcommand{\ds} {\displaystyle}

\newcommand{\RR} {\mathbb{R}}
\newcommand{\CC} {\mathbb{C}}
\newcommand{\DD} {\mathbb{D}}

\newcommand{\FF} {\mathbb{F}}
\newcommand{\LL} {\mathbb{L}}
\newcommand{\TT} {\mathbb{T}}
\newcommand{\ZZ} {\mathbb{Z}}
\newcommand{\bC} {\mathbb{C}}

\newcommand{\ptq}{p\times q}
\newcommand{\qtp}{q\times p}

\newcommand{\cB} {\mathcal{B}}
\newcommand{\cC} {\mathcal{C}}

\newcommand{\cS} {\mathcal{S}}

\newcommand{\cH} {\mathcal{H}}

\newcommand{\cM} {\mathcal{M}}
\newcommand{\cN} {\mathcal{N}}

\newcommand{\cX}{\mathcal{X}}
\newcommand{\cY}{\mathcal{Y}}

\newcommand{\cW}{\mathcal{W}}
\newcommand{\cZ}{\mathcal{Z}}

\newcommand{\mtm} {m\times m}
\newcommand{\ptp} {p\times p}

\newcommand{\gE} {\mathfrak{E}}

\newcommand{\gF} {\mathfrak{F}}
\newcommand{\gee}{\mathfrak{e}}

\newcommand{\gu}{\mathfrak{u}}

\newcommand{\gv}{\mathfrak{v}}
\newcommand{\gw}{\mathfrak{w}}

\newcommand{\gV}{\mathfrak{V}}

\newcommand{\gz}{\mathfrak{z}}

\newcommand{\dsp}{\displaystyle}
\def\wt{\widetilde}
\def\wh{\widehat}
\def\ol{\overline}
\def\CR{\color{red} }
\def\CB{\color{blue} }

\title[ ]{Toeplitz, Hankel, de Branges and two truncated matrix moment problems}

\author[Dhara]{Kousik Dhara}
\address{Kousik Dhara, Department of Mathematics, Indian Institute of Technology Jammu, Jammu and Kashmir 181221, India}
\email{kousik.dhara@iitjammu.ac.in}

\author[Dym]{Harry Dym$^\dagger$}
\address{Harry Dym, Department of Mathematics, Weizmann Institute of Science, Rehovot 7610001, Israel}
\email{harry.dym@weizmann.ac.il}


\begin{document}
	
\dedicatory{To our distinguished colleague Damir Z. Arov on the occasion of his ninetieth birthday}	
\thanks{$\dagger$ Professor Harry Dym sadly passed away on July 18, 2024, shortly after the completion of
this article.
}	
	
	\begin{abstract}
This paper deals with
(1) the truncated matrix Hamburger moment problem from the point of view of reproducing kernel Hilbert spaces of vector valued entire functions of the kind introduced and extensively studied by Louis de Branges and (2) the truncated matrix trigonometric moment problem viewed through an analogous class of spaces that are formulated with respect to the open unit disc rather than the open upper half-plane. In this approach projections are computed via appropriately chosen reproducing kernels instead of orthogonal bases. This approach eases the bookkeeping and leads to pleasing formulas.	
	\end{abstract}
	\maketitle

	\tableofcontents

	\section{Introduction}
	 This paper deals with 
  the truncated matrix Hamburger moment problem 
	 and the truncated matrix trigonometric moment problem from what seems to be a new point of view: the intimate connection of the first  of these two problems with a special class of RKHS's (reproducing kernel Hilbert spaces) that was introduced and extensively studied by L. de Branges and the second with an analogous class of RKHS's 
	 	 that is formulated with respect to the open unit disc $\DD$ rather than the open upper half-plane $\CC_+$. To explain this connection in more detail, let 
	 $\CC^{\ptq}$ denote the set of $p\times q$ matrices with complex entries, $\CC^p=\CC^{p\times 1}$ and $\CC=\CC^1$ and let
	\begin{equation}
		\label{eq:aug24a23}
		F(\lambda)=\begin{bmatrix}I_p&\lambda I_p&\cdots&\lambda^nI_p\end{bmatrix}\quad\textrm{for}\quad 
		\lambda\in\CC,\ p\ge 1, \ n\ge 1.
	\end{equation}
	Then the space
	\begin{equation}
		\label{eq:aug24b23}
		{\cH}_G=\{Fu:\,u\in\CC^{(n+1)p}\}
	\end{equation}
	endowed with the inner product
	\begin{equation}
		\label{eq:aug24c23}
		\langle Fu, Fv\rangle_{{\cH}_G}=v^*Gu
	\end{equation}
	based on the positive definite matrix
	\begin{equation}
		\label{eq:aug24d23}
		G=\begin{bmatrix}g_{00}&\cdots&g_{0n}\\ \vdots&&\vdots\\
			g_{n0}&\cdots&g_{nn}\end{bmatrix}\quad\textrm{with}\ g_{ij}\in\CC^{\ptp}\ \textrm{for}\ i,j=0,\ldots,n
	\end{equation}
	is a Hilbert space of $p\times 1$ vector valued polynomials. It is readily checked that if 
	\begin{equation}
		\label{eq:aug24e23}
		\Gamma=G^{-1}\quad\textrm{and}\quad K_\omega(\lambda)=F(\lambda)\Gamma F(\omega)^*,
	\end{equation}
	then, for every choice of $x\in\CC^p$, $u\in\CC^{(n+1)p}$  and $\omega\in\CC$, 
	\begin{equation}
	\label{eq:aug24f23}
	K_\omega x\in{\cH}_G\quad\textrm{and}\quad 
	\langle Fu,K_\omega x\rangle_{{\cH}_G}=x^*F(\omega)u.
\end{equation}
Thus, ${\cH}_G$ is an RKHS  with RK (reproducing kernel) $K_\omega(\lambda)$.

In the sequel, we shall be especially interested in the special class of RKHS's for which the RK can be expressed in terms of a  
$p\times 2p$ mvf (matrix valued function) $\gE(\lambda)=\begin{bmatrix}E_-(\lambda)&E_+(\lambda)\end{bmatrix}$ with blocks $E_\pm(\lambda)$ of size $\ptp$ by the formula
\begin{equation}
	\label{eq:aug24g23}
	K_\omega(\lambda)=\left\{\begin{array}{cc}
		\frac{\dsp E_+(\lambda)E_+(\omega)^*-E_-(\lambda)E_-(\omega)^*}{\dsp -2\pi i(\lambda-\ol{\omega})}&\quad\textrm{if $\lambda\ne\ol{\omega}$}\\ 
		\frac{\dsp E_+^\prime(\ol{\omega})E_+(\omega)^*-E_-^\prime(\ol{\omega})E_-(\omega)^*}{\dsp -2\pi i
		}&\quad\textrm{if $\lambda=\ol{\omega}$}
	\end{array}\right.
\end{equation}
and shall refer to such RKHS's as de Branges spaces ${\cB}(\gE)$ with respect to $\CC_+$, the open upper half-plane. 
This class of spaces  of entire vvf's was introduced by L. de Branges in \cite{br63} and  extensively studied by him in a series of  papers in the late fifties and sixties that culminated in the monograph \cite{br68a} for the case $p=1$.  Generalizations to the case $p>1$ were considered  in \cite{brr66} and \cite{br68b} and to wider classes of spaces in numerous other publications.

We shall also consider an analogous class of spaces with RK's of the form 
\begin{equation}
	\label{eq:aug24h23}
	K_\omega(\lambda)=\left\{\begin{array}{cc}
		\frac{\dsp E_+(\lambda)E_+(\omega)^*-E_-(\lambda)E_-(\omega)^*}{\dsp 1-\lambda \ol{\omega}}&\quad\textrm{if}\ \lambda\ol{\omega}\ne 1\\
	-\lambda E_+^\prime(\lambda)E_+(\omega)^*+\lambda E_-^\prime(\lambda)E_-(\omega)^*&\quad\textrm{if}\ \lambda\ol{\omega}=1,
	\end{array}\right.
\end{equation}
and shall refer to this class as de Branges spaces with respect to $\DD$, the open unit disk. 

This paper rests on the observation that:
\begin{align}
\label{eq:dec18a23}
&{\cH}_G\ \textrm{is a de Branges space with respect to $\CC_+$ $\iff$ $G$ is a block Hankel matrix.}\\
\label{eq:dec18b23}
&{\cH}_G\ \textrm{is a de Branges space with respect to $\DD$ $\iff$ $G$ is a block Toeplitz matrix.}
\end{align}

The fact that ${\cH}_G$ can be identified as a de Branges space when $G$ is block Toeplitz or block Hankel has been observed before, see e.g., \cite{d88}, \cite{d89a}, \cite{bd96}. The converse implication   $\Longrightarrow$  seems to 
have been overlooked (apart from a brief discussion of the scalar case in \cite{d23b}).

The verification of these statements   depends on a characterization of de Branges spaces of vector valued entire functions with respect to $\CC_+$ that originates with de Branges for scalar valued entire functions in \cite{br68a} and was subsequently observed to hold for $p>1$ in \cite{ds17}, and on an analogous characterization of de Branges spaces with respect to $\DD$ that was presented in \cite{d23a}. 
	The proof is constructive in the sense that it presents formulas for $E_\pm$ in both settings.

There is a natural connection between the two classes of de Branges spaces under consideration and two truncated matrix moment problems. Recall that  nondecreasing mvf $\sigma(t)$ on $\RR$ is said to be a solution of the truncated matrix Hamburger moment problem with matrix moments $h_0,\ldots,h_{2n}\in\CC^{\ptp}$ if
\begin{equation}
\label{eq:dec7b23}
\int_{-\infty}^\infty \mu^k d\sigma(\mu)=h_k \quad\textrm{for} \ k=0,\ldots, 2n.
\end{equation}
Analogously, a  nondecreasing mvf $\sigma(t)$ on $[0,2\pi]$ is said to be a solution of the truncated matrix trigonometric  moment problem with matrix moments $h_{-n},\ldots,h_{n}\in\CC^{\ptp}$ if
\begin{equation}
\label{eq:dec7c23}
\frac{1}{2\pi}\int_0^{2\pi} e^{-ikt} d\sigma(t)=h_k \quad\textrm{for} \ k=-n,\ldots, n.
\end{equation}
 It is well known (see e.g., Corollary 3 to Theorem 1 and Theorem 4 in \cite{an70}) that: 
 \[
 \begin{split}
&\textrm{A necessary  condition for the existence of a solution to \eqref{eq:dec7b23} (resp., \eqref{eq:dec7c23}) is that}\\ &\textrm{the block Hankel matrix 
$G$ with entries $g_{ij}=h_{i+j}$ (resp., the block Toeplitz matrix $G$ with entries}\\ &\textrm{$g_{ij}=h_{i-j}$) for $i,j=0,\ldots,n$ be positive semidefinite.}
\end{split}
\] 
Thus, to obtain a solution, we must choose $h_k=h_k^*$ for $k=0,\ldots,2n$ in \eqref{eq:dec7b23} and $h_k=h_{-k}^*$ for $k=-n,\ldots,n$ in \eqref{eq:dec7c23}. 
 
In this article we restrict attention to positive definite matrices $G$. This ensures that 
 there exists an absolutely continuous solution $\sigma(t)$ with density $\sigma^\prime =(E_+E_+^*)^{-1}$ based on the mvf $E_+$ that intervenes in formulas  \eqref{eq:aug24g23} and \eqref{eq:aug24h23}, respectively.   This is an easy byproduct of the analysis in Sections \ref{sec:toeplitz} and \ref{sec:hankel}, which are devoted  to the justification of \eqref{eq:dec18a23} and \eqref{eq:dec18b23}. These two sections are preceded by 
 Section \ref{sec:prelim}, which presents a brief review of some basic facts about RKHS's and some relevant function spaces for the convenience of the reader.  
 
   Section \ref{sec:dbspacecharacterizations} is a short detour to  briefly explain the basic idea that underlies the characterization of de Branges spaces that originates in \cite{br68a} and deserves to be known better. Therein it is also shown that  the isometries that serve to characterize the spaces under consideration are equivalent to  important matrix identities of the Gohberg-Heinig, Gohberg-Semencul type. 
  
Sections \ref{sec:solutionstrig} and \ref{sec:solutionshamburger} are devoted to obtaining the set of all solutions to the  two moment problems under consideration. The starting point in both of these sections is  an identity that expresses the set of all solutions of the moment problem in terms of the Arov-Grossman functional model of the set of minimal isometric extensions of an appropriately chosen  isometric operator. In this we follow in the 
 footsteps of \cite{mm03} and  \cite{mm05} which treated the scalar case and the matrix case of the truncated matrix trigonometric moment problem, respectively and \cite{za12} which treated the truncated matrix Hamburger moment problem. The trigonometric problem was also treated by operator theoretic methods in \cite{za11} and \cite{za13}. However, we compute the requisite projections via reproducing kernels, rather than orthogonal bases.  (See e.g., \eqref{eq:feb26b24}.) This eases the bookkeeping and yields simple formulas for the entries $\Theta_{ij}$ in the linear fractional transformation 
 $$
T_\Theta[S]=(\Theta_{11}S+\Theta_{12})(\Theta_{21}S+\Theta_{22})^{-1}
$$
that is used to parameterize the solutions of both of the moment problems in terms of mvf"s 
$S$ in the Schur class ${\cS}^{\ptp}(\Omega_+)$ of $\ptp$ mvf's that are holomorphic and contractive in $\Omega_+$, $\Omega_+$ is either $\DD$ or $\CC_+$. The blocks $\Theta_{ij}$ are computed from recipes based on $E_\pm$; 
see Theorems \ref{thm:may7a24},  \ref{thm:jun5a24} and Section \ref{sec:epilogue}, which is devoted to summarizing remarks and a brief discussion of supplementary information.

\textbf{Literature:} There is an extensive literature on truncated moment problems. One of the best general references is still the classic monograph \cite{ak65} by Akhiezer, supplemented by \cite{bro71} to help convert conclusions for scalar problems to matrix problems;   \cite{za12} has a good discussion of literature.

 This article, however, was largely motivated by a number of the formulas in the treatment  of the Hamburger moment problem in\cite{ga92} and in the treatment of the trigonometric moment problem in \cite{mm03}, which partially overlapped 
formulas that arise in the theory of  scalar 
de Branges  spaces. This lead us to think that these spaces would be a natural tool for the analysis of such moment problems 
and in fact for matrix versions of these problems, which is what this article is about.

\textbf{Notation:} The notation is mostly standard and self-evident:
The symbols $\CC$, $\CC^{\ptq}$ and $\CC^p$ will be used to denote the complex numbers, the set of $\ptq$ matrices with complex entries and the set of $p\times 1$ vectors with complex entries, respectively; $\CC_-$ denotes the open lower half-plane.  The symbol $\RR$ denotes the real numbers and $\TT=\{\lambda\in\CC:\,\vert\lambda\vert=1\}$;   if $A\in\CC^{\ptq}$ then $A^T$ and $A^*$ denote the  transpose and the Hermitian transpose of $A$ (alias the adjoint of $A$ with respect to the standard inner product $\langle{u},{v}\rangle_{st}={v}^*{u}$ in $\CC^p$), respectively. 
The letter $J$ is used to denote a {\bf signature matrix}, i.e., a matrix that is both selfadjoint and unitary with respect to the standard inner product. The particular signature matrices  
\begin{equation}
\label{eq:30d23}
	J_p=\begin{bmatrix}
		0 & -I_p\\
		-I_p & 0
	\end{bmatrix}, \quad
	j_p=\begin{bmatrix}
		I_p & 0 \\
		0 & -I_p
	\end{bmatrix} ,\quad {\gV}=\frac{1}{\sqrt{2}}\begin{bmatrix}-I_p&I_p\\I_p&I_p\end{bmatrix},
\end{equation}
which are connected by the formulas 
\begin{equation}
\label{eq:30e23}
{\gV}J_p{\gV}=j_p\quad\textrm{and}\quad {\gV}j_p{\gV}=J_p
\end{equation}
will be used.

Throughout this paper 
$$
p\quad\textrm{and}\quad n\quad \textrm{will be fixed positive integers and}\quad m=(n+1)p,
$$
\[
F(\lambda)=\begin{bmatrix}I_p&\lambda I_p&\cdots&\lambda^nI_p\end{bmatrix}
\quad\textrm{and}\quad\wh{F}(\lambda)=\lambda F(\lambda),
\]
$G\in\CC^{\mtm}$  will be a positive definite block matrix with entries $g_{ij}\in\CC^{\ptp}$ for $i,j=0,\ldots,n$,

\begin{align*}
	G_{[j,k]}=\begin{bmatrix}
		& g_{jj} & \cdots & g_{jk} \\
		&\vdots &         & \vdots \\
		&g_{kj} & \cdots & g_{kk} 
	\end{bmatrix}\quad 
	\text{for}\  0\le j\leq k\leq n
\end{align*}
and
$$
G^{-1}=\Gamma=\begin{bmatrix}\gamma_{00}^{(n)}&\gamma_{01}^{(n)}&\cdots&\gamma_{0n}^{(n)}\\
\vdots& & &\vdots\\
\gamma_{n0}^{(n)}&\gamma_{n1}^{(n)}&\cdots&\gamma_{nn}^{(n)} \end{bmatrix}.
$$

Since polynomials of degree $\le n$ include polynomials of degree $0$ it is convenient to set 
$$
{\gee_k}=\begin{bmatrix}0&\cdots&0&I_p&0&\cdots&0\end{bmatrix}^T\quad \textrm{for}\ k=0,\ldots,n 
$$
where the $I_p$ appears as the $k+1$ entry in the $m\times p$ matrix  ${\gee}_k$, and to let $\{\eta_1,\ldots,\eta_p\}$ denote the standard basis for $\CC^p$. Then the set of $m=(n+1)p$ vectors 
$$
\{{\gee}_j\eta_k\},\  j=0,\ldots,n,\  k=1,\ldots,p,
$$
 is an orthonormal basis for $\CC^m$ with respect to the standard inner product. We shall set 
 $$
{\gu}_j=\Gamma {\gee}_j (\gamma_{jj}^{(n)})^{-1/2}, \quad \textrm{for}\ j=0,\ldots,n
$$
 and shall use the symbol $A$ to  denote 
the $(n+1)p\times(n+1)p$ matrix of the form: 
$$
A=\sum_{k=0}^{n-1}{\gee}_j{\gee}_{j+1}^T
$$
Thus, for example, if $n=3$, then
$$
A=\sum_{j=0}^3{\gee}_j{\gee}_{j+1}^T=\begin{bmatrix}
	0&I_p&0&0\\ 0&0&I_p&0\\ 0&0&0&I_p\\ 0&0&0&0\end{bmatrix}.
$$

Some additional notation will be introduced in the next section.

\section{Preliminaries}
\label{sec:prelim}
In this section we shall review some needed facts about RKHS's and related  function spaces.
\subsection{RKHS's} 
It is well known (and easily verified) that 
if ${\cH}$ is a RKHS of entire $p\times 1$ vvf's with RK's $K_\omega(\lambda)$ and $L_\omega(\lambda)$, then:
\begin{equation}
\label{eq:nov17a23}
\Vert K_\omega x-L_\omega x\Vert_{\cH}=0\quad\textrm{for every choice of $\omega\in\CC$ and $x\in\CC^p$}
\end{equation}
and
\begin{equation}
\label{eq:nov17b23}
\sum_{i,j=1}^n x_i^* K_{\omega_j}(\omega_i)x_i=\Vert \sum_{i=1}^n K_{\omega_j} x_j\Vert_{\cH}^2\ge 0
\end{equation}
for every choice of $\omega_1,\ldots,\omega_n\in\CC$ and $x_1,\ldots,x_n\in\CC^p$.
 Thus,
\begin{enumerate}
\item[\rm(1)]  A RKHS ${\cH}$ has exactly one RK.
\item[\rm(2)] 
 $K_\omega(\lambda)$ is a positive kernel in the sense that the double sum  in \eqref{eq:nov17b23} is always $\ge 0$.
\item[\rm(3)] $K_\alpha(\beta)=K_\beta(\alpha)^*$. (This follows from (2) with $n=2$.)
\end{enumerate}

The subspace
\begin{equation}
		\label{eq:nov17c23}
				{\cH}_\omega=\{f\in{\cH}:\ f(\omega)=0\}
	\end{equation}
and the operator
\begin{equation}
\label{eq:nov15e23}
(R_\alpha f)(\lambda)=\left\{\begin{array}{cc}\frac{\ds f(\lambda)-f(\alpha)}{\ds \lambda-\alpha} &\quad\textrm{if}\ \lambda\ne \alpha\\ f^\prime(\alpha)&\quad\textrm{if}\ \lambda=\alpha
\end{array}\right.
\end{equation}
will play an important role in the sequel. 

 In terms of the notation $F(\lambda)$, ${\gee}_j$, $m$ and $A$ that was introduced in the last section, it is readily checked that 
\begin{equation}
\label{eq:nov17e23}
F(\lambda)={\gee}_0^*(I_m-\lambda A)^{-1}
\end{equation}
and hence that the space 
${\cH}_G$ of vector valued polynomials with inner product \eqref{eq:aug24c23} based on a positive definite matrix $G\in\CC^{\mtm}$ can be described as 
\begin{equation}
\label{eq:feb25a24}
{\cH}_G=\{F u:\ u\in\CC^m\} \quad\textrm{and}\quad \langle Fu,Fv\rangle_{{\cH}_G}=v^*Gu \quad\text{for every}\   u,v\in\CC^m.
\end{equation}
The main properties of ${\cH}_G$ are summarized in the next theorem. They are partially expressed in terms of the {\bf notation}
\begin{equation}
\label{eq:feb26a24}
{\gz}_\alpha=\Gamma F(\alpha)^*(F(\alpha)\Gamma F(\alpha)^*)^{-1/2}  \quad \textrm{and}\quad      N_\alpha=\Gamma-{\gz}_\alpha {\gz}_\alpha^*.
\end{equation}
\begin{thm}
\label{thm:feb25a24}
If $G$ is a positive definite matrix, then ${\cH}_G$ is a RKHS with RK
\begin{equation}
\label{eq:feb25b24}
K_\omega(\lambda)=F(\lambda)\Gamma F(\omega)^* \quad\text{for every choice of}\ \lambda,\omega\in\CC.
\end{equation}
Moreover, 
\begin{enumerate}
\item[\rm(1)] $K_\alpha(\alpha)=F(\alpha)\Gamma F(\alpha)^*\succ O$ for every point $\alpha\in\CC$.
\item[\rm(2)] ${\cH}_G$ is $R_\alpha$ invariant for every point $\alpha\in\CC$ and 
\begin{equation}
\label{eq:nov17f23}
(R_\alpha F)(\lambda)= F(\lambda)A(I_m-\alpha A)^{-1}.
\end{equation}
\item[\rm(3)]  $({\cH}_G)_\alpha=\{f\in{\cH}_G:\ f(\alpha)=0\}$ is a RKHS with RK 
\begin{equation}
\label{eq:nov22e23}
K_\omega^{(\alpha)}(\lambda)=
F(\lambda)N_\alpha F(\omega)^*.
\end{equation}
\item[\rm(4)] ${\cH}_G=({\cH}_G)_\alpha\oplus \{K_\alpha u: \ u\in\CC^p\}$.
\item[\rm(5)] $({\cH}_G)_\alpha=\{FN_\alpha u:\,u\in\CC^m\}$.
\item[\rm(6)] The orthogonal projections $\Pi_\alpha$ of ${\cH}_G$ onto $({\cH}_G)_\alpha$ and $\Pi_\alpha^\perp$ of 
${\cH}_G$ onto ${\cH}_G\ominus ({\cH}_G)_\alpha$ are given by the formulas
\begin{equation}
\label{eq:feb26b24} 
\Pi_\alpha Fu=FN_\alpha Gu\quad and \quad\Pi_\alpha^\perp Fu=F{\gz}_\alpha {\gz}_\alpha^*Gu. 
\end{equation}
\item[\rm(7)] $\{(I_m-\alpha A)^{-1}N_\alpha u:\,u\in\CC^m\}=\{A^*Au:\,u\in\CC^m\}=\{A^*u:\,u\in\CC^m\}$.
\item[\rm(8)] $\det\,K_\omega(\lambda)$ is a polynomial of degree $\le np$ in $\lambda$ and 
\begin{equation}
\label{eq:mar4a24}
K_\omega(\lambda)^{-1}=K_\omega(\omega)^{-1}F(\omega){\gz}_\omega\{F(\lambda){\gz}_\omega\}^{-1}
\end{equation}
when the indicated inverse exists.
\end{enumerate}
\end{thm}

\begin{proof}
In view of  \eqref{eq:feb25a24}, 
$$
\langle Fu, F\Gamma F(\omega)^*v\rangle_{{\cH}_G}=v^*F(\omega)\Gamma Gu=v^*F(\omega)u
$$
for every choice of $u,v\in\CC^m$ and $\omega\in\CC$. Therefore, \eqref{eq:feb25b24} holds. Moreover, $K_\alpha(\alpha)x=0\implies F(\alpha)^*x=0\implies x=0$, since $\Gamma$ is positive definite and $\textup{rank}\,F(\alpha)=p$. Therefore, (1) holds.

Item (2) follows from the fact that 
\[
\begin{split}
(I_m-\lambda A)^{-1}-(I_m-\alpha A)^{-1}&=(I_m-\lambda A)^{-1}\{(I_m-\alpha A)-(I_m-\lambda A)\}(I_m-\alpha A)^{-1}\\
&=(\lambda-\alpha)(I_m-\lambda A)^{-1}A(I_m-\alpha A)^{-1}.
\end{split}
\]
Moreover, since the kernel 
\[
K_\omega^{(\alpha)}(\lambda)=K_\omega(\lambda)-K_\alpha(\lambda)K_\alpha(\alpha)^{-1}K_\omega(\alpha)=F(\lambda)N_\alpha F(\omega)^*
\]
meets the conditions $K_\omega^{(\alpha)}u\in({\cH}_G)_\alpha$ and 
$$
\langle f, K_\omega^{(\alpha)}v\rangle_{{\cH}_G}=v^*f(\omega)\quad\textrm{for every}\ f\in({\cH}_G)_\alpha,  
$$
it is the one and only  RK for $({\cH}_G)_\alpha$.  Therefore (3) holds. 

Next, (4) is self-evident and  it is also clear that $\{FN_\alpha u:\,u\in\CC^m\}\subseteq ({\cH}_G)_\alpha$, since
\begin{equation}
\label{eq:feb27a24}
F(\alpha)N_\alpha=F(\alpha)(\Gamma-{\gz}_\alpha{\gz}_\alpha^*)=F(\alpha)(\Gamma-\Gamma F(\alpha)^*(F(\alpha)\Gamma F(\alpha)^*)^{-1}
F(\alpha)\Gamma=0. 
\end{equation}
On the other hand, if $f=Fu$ belongs to $({\cH}_G)_\alpha$ and $\langle f,FN_\alpha v\rangle_{{\cH}_G}=0$ for every $v\in\CC^m$, then $F(\alpha)u=0$ and 
$$
0=(\Gamma-{\gz}_\alpha{\gz}_\alpha^*)Gu=u-\Gamma F(\alpha)^*(F(\alpha)\Gamma F(\alpha)^*)^{-1}F(\alpha)u=u.
$$
Therefore, equality prevails in the preceding inclusion, i.e., (5) holds. 

To verify (6), observe that
\[
\begin{split}
v^* (\Pi_\alpha Fu)(\omega)&=\langle 
\Pi_\alpha Fu, K_\omega^{(\alpha)}v\rangle_{{\cH}_G}=\langle Fu, K_\omega^{(\alpha)}v\rangle_{{\cH}_G}\\
&=\langle Fu, FN_\alpha F(\omega)^*v\rangle_{{\cH}_G}=v^*F(\omega)N_\alpha Gu
\end{split}
\]
for every $v\in\CC^m$. Therefore, 
 the first formula in    \eqref{eq:feb26b24} holds; the second follows from the first, since 
$\Pi_\alpha^\perp=I-\Pi_\alpha$. 

Next, in view of \eqref{eq:feb27a24} and the identity $I_m={\gee}_0{\gee}_0^*+A^*A$,
$$
(I_m-\alpha A)^{-1}N_\alpha=({\gee}_0{\gee}_0^*+A^*A)(I_m-\alpha A)^{-1}N_\alpha=A^*A(I_m-\alpha A)^{-1}N_\alpha.
$$
Therefore,  $\{(I_m-\alpha A)^{-1}N_\alpha u:\,u\in\CC^m\}\subseteq\{A^*Au:\,u\in\CC^m\}$.
 Moreover, as 
\[
(I_m-\alpha A)^{-1}N_\alpha u=0\iff N_\alpha u=0\iff  (I_m-F(\alpha)^*K_\alpha(\alpha)^{-1}F(\alpha)\Gamma )u=0
\]
it follows that the dimension of the null space $(I_m-\alpha A)^{-1}N_\alpha$ is $\le p$ and hence that
$$
\textup{rank}\,(I_m-\alpha A)^{-1}N_\alpha \ge np=\textup{rank}\,A^*A=\textup{rank}\,A^*.
$$
Thus, (7) holds.

To verify (8),
\[
\begin{split}
K_\omega(\lambda)&=\{F(\lambda)-F(\omega)\}\Gamma F(\omega)^*+K_\omega(\omega)\\
&=(\lambda-\omega)F(\lambda)A(I_m-\omega A)^{-1}\Gamma F(\omega)^*+K_\omega(\omega)\\
&=\{I_p+(\lambda-\omega)F(\lambda)A(I_m-\omega A)^{-1}\Gamma F(\omega)^*K_\omega(\omega)^{-1}\} K_\omega(\omega)\\
&=\{I_p-(\omega-\lambda)F(\omega)A(I_m-\lambda A)^{-1}\Gamma F(\omega)^*K_\omega(\omega)^{-1}\} K_\omega(\omega)
\end{split}
\]
The term in curly brackets is of the form $I_p-BC$ with $B=(\omega-\lambda)F(\omega)A(I_m-\lambda A)^{-1}$ 
and $C=\Gamma F(\omega)^*K_\omega(\omega)^{-1}={\gz}_\omega K_\omega(\omega)^{-1/2}$. Therefore, 
since $\det(I_p-BC)\not\equiv 0$ and ${\gz}_\omega^*G{\gz}_\omega=I_p$, 
\[
\begin{split}
(I_p-BC)^{-1}&=I_p+B(I_m-CB)^{-1}C\\ &=I_p+B\{I_m-(\omega-\lambda){\gz}_\omega{\gz}_\omega^*GA(I_m-\lambda A)^{-1}\}^{-1}{\gz}_\omega K_\omega(\omega)^{-1/2}\\
&=I_p+B{\gz}_\omega\{I_p-(\omega-\lambda){\gz}_\omega^*GA(I_m-\lambda A)^{-1}{\gz}_\omega\}^{-1}K_\omega(\omega)^{-1/2}\\
&=I_p+B{\gz}_\omega\{{\gz}_\omega^*G[I_m-(\omega-\lambda)A(I_m-\lambda A)^{-1}]{\gz}_\omega\}^{-1}K_\omega(\omega)^{-1/2}\\
&=I_p+B{\gz}_\omega\{{\gz}_\omega^*G[(I_m-\omega A)(I_m-\lambda A)^{-1}]{\gz}_\omega\}^{-1}K_\omega(\omega)^{-1/2}\\
&=I_p+B{\gz}_\omega\{K_\omega(\omega)^{-1/2}F(\lambda){\gz}_\omega\}^{-1}K_\omega(\omega)^{-1/2}\\
&=I_p+B{\gz}_\omega\{ F(\lambda){\gz}_\omega\}^{-1}\\
&=F(\omega){\gz}_\omega \{F(\lambda){\gz}_\omega\}^{-1}
\end{split}
\]
for all points $\lambda\in\CC$ for which the indicated inverse exists. 
Consequently, \eqref{eq:mar4a24} holds.
 \end{proof}
 
 It is also convenient to assemble a number of properties of the matrix $N_\alpha$ for future use.
 \begin{thm}
 \label{thm:may14a24}
 If $G$ is a positive definite matrix, then
 \begin{enumerate}
 \item[\rm(1)] ${\gz}_\alpha^*G{\gz}_\alpha=I_p$ and hence $\textup{rank}\,{\gz}_\alpha=p$.
 \item[\rm(2)] $N_\alpha G N_\alpha=N_\alpha$ (and hence $N_\alpha$ is positive semi-definite) and $N_\alpha G{\gz}_\alpha=0$.
 \item[\rm(3)] $F(\alpha) N_\alpha=0$.
 \item[\rm(4)] $\textup{rank}\,N_\alpha=m-p$.
 \item[\rm(5)] $\Vert G^{1/2}N_\alpha G^{1/2}\Vert= 1$.
 \end{enumerate}
 \end{thm}
 
 \begin{proof} The first three assertions are straightforward computations that are left to the reader. The fourth is obtained by using Schur complement formulas to check that the three matrices
 \[
 \begin{bmatrix}
 \Gamma&\Gamma F(\alpha)^*\\ 
 F(\alpha)\Gamma&K_\alpha(\alpha)\end{bmatrix},\quad \begin{bmatrix}N_\alpha&O\\ O&K_\alpha(\alpha)\end{bmatrix},\quad \begin{bmatrix} \Gamma&O\\O&O\end{bmatrix}
 \]
 are congruent and hence must have the same rank.
 
To verify (5), observe that
\[
{u}^*GN_\alpha GN_\alpha G{u}=\Vert FN_\alpha Gu\Vert_{{\cH}_G}^2=\Vert \Pi_\alpha Fu\Vert_{{\cH}_G}^2\le \Vert Fu\Vert_{{\cH}_G}^2
=u^*Gu
\]
for every $u\in\CC^m$. 
Thus, as $(N_\alpha G)^2=N_\alpha G$, $GN_\alpha G\preceq G$ and  $\Vert G^{1/2}N_\alpha G^{1/2}\Vert\le 1$. 
On the other hand, the identity $G^{1/2}N_\alpha G^{1/2}G^{1/2}N_\alpha=G^{1/2}N_\alpha$ exhibits the nonzero columns 
of  $G^{1/2}N_\alpha$ as eigenvectors of $G^{1/2}N_\alpha G^{1/2}$ corresponding to the eigenvalue $1$. Therefore, 
$\Vert G^{1/2}N_\alpha G^{1/2}\Vert\ge 1$ and hence (5) holds.


  \end{proof}
For additional information on RKHS's of entire vvf's see e.g., [ArD18]. 

\subsection{Function spaces} 
The symbol $\Omega_+$ stands for either $\CC_+$ or $\DD$, and correspondingly $\Omega_0$ stands for either $\RR$ or $\TT$ (alias the boundary of $\Omega_+$). 
We shall use the symbols $L_1(\Omega_0)$, $L_2(\Omega_0)$ and $L_\infty(\Omega_0)$ to denote the spaces of measurable functions on $\Omega_0$ that are summable, square summable and bounded, respectively; $L_r^{\ptq}(\Omega_0)$ will denote the set of $\ptq$ mvf's with norms that belong to $L_r(\Omega_0)$.

The Hardy spaces  $H_2^p(\DD)$ (resp., $H_2^p(\CC_+)$) of $p\times 1$ vvf's $f$ that are holomorphic in $\DD$ (resp., $\CC_+$) that meet the condition
$$
\sup_{0\le r<1}\int_0^{2\pi}f(re^{it})^*f(re^{it})dt<\infty\quad\textrm{(resp., $\sup_{\nu>0}\int_{-\infty}^\infty f(\mu+i\nu)^*f(\mu+i\nu)d\mu<\infty$)}
$$
are  RKHS's of vvf's that are holomorphic in $\Omega_+$ with RK 
\begin{equation}
\label{eq:nov21a23}
K_\omega(\lambda)=\frac{I_p}{\rho_\omega(\lambda)},\quad \textrm{where}\ \rho_\omega(\lambda)=\left\{\begin{array}{cc}
1-\lambda\ol{\omega}&\ \textrm{if}\ \Omega_+=\DD\\
-2\pi i(\lambda-\ol{\omega})&\quad\textrm{if}\ \Omega_+=\CC_+\end{array}\right. 
\end{equation}
and inner product $\langle f,g\rangle_{\rm st}$ that is defined in terms of the non tangential boundary values by the formulas
\begin{equation}
\label{eq:nov21b23}
\langle f,g\rangle_{\rm st}=\left\{\begin{array}{ll}\frac{1}{2\pi}\int_0^{2\pi}g(e^{it})^*f(e^{it})dt &\quad \textrm{if}\ \Omega_+=\DD\\ \\
\int_{-\infty}^\infty g(\mu)^*f(\mu)d\mu&\quad \textrm{if}\ \Omega_+=\CC_+.\end{array}\right.
%
\end{equation}
Correspondingly, the orthogonal projection $\Pi_+f$ of a vvf $f\in L_2^p(\TT)$ (resp., $f\in L_2^p(\RR)$) 
is obtained from the observation that 
\begin{equation*}
\label{eq:nov21c23}
\xi^*(\Pi_+f)(\omega)=\langle \Pi_+f,(\rho_\omega)^{-1}\xi\rangle_{\rm st}=\langle f, \Pi_+(\rho_\omega)^{-1}\xi\rangle_{\rm st}
=\langle f,(\rho_\omega)^{-1}\xi\rangle_{\rm st}
\end{equation*}
for every choice of $\xi\in\CC^p$ and $\omega\in\Omega_+$ and hence that 
\begin{equation}
\label{eq:nov21c23}
\begin{split}
(\Pi_+f)(\omega)&=\frac{1}{2\pi}\int_0^{2\pi}\frac{f(e^{it})}{1-e^{-it}\omega}dt\quad\textrm{if $\omega\in\Omega_+$ and $\Omega_+=\DD$}\quad \textrm{and}\\
(\Pi_+f)(\omega)&=\frac{1}{2\pi i}\int_{-\infty}^{\infty}\frac{f(\mu)}{\mu-\omega}d\mu\quad\textrm{if $\omega\in\Omega_+$ and $\Omega_+=\CC_+$}.
\end{split}
\end{equation}
When $f=\Psi u$ for some $\ptp$ mvf $\Psi$ with square summable entries and $u\in\CC^p$ we shall often write
\begin{equation}
\label{eq:nov30e23}
\Pi_+\Psi \quad\textrm{understanding that the projection is applied column by column, instead of}\ 
\Pi_+\Psi u,
\end{equation}
in order to minimize clutter.

The following subspaces of the class ${\cX}^{\ptp}(\Omega_+)$ of $\ptp$ mvf's that are holomorphic in $\Omega_+$ will come into play:
\[
\begin{split}
H_\infty^{\ptp}(\Omega_+)&=\{f\in {\cX}^{\ptp}(\Omega_+):\, \sup_{\omega\in\Omega_+} \Vert f(\omega)\Vert<\infty\}\\
{\cS}^{\ptp}(\Omega_+)&=\{f\in {\cX}^{\ptp}(\Omega_+):\,f(\omega)^*f(\omega)\preceq I_p\quad\textrm{for}\ \omega\in\Omega_+\}\quad\textrm{(The Schur class)}\\
{\cS}_{\rm in}^{\ptp}(\Omega_+)&=\{f\in{\cS}^{\ptp}(\Omega_+):\ f(\omega)^*f(\omega)=I_p\quad\textrm{for  
almost all points $\omega\in\Omega_0$}\}\\
{\cC}^{\ptp}(\Omega_+)&=\{f\in {\cX}^{\ptp}(\Omega_+):\, f(\omega)+f(\omega)^*\succeq 0\quad\textrm{for}\ \omega\in\Omega_+\}\quad\textrm{(The Carath\'{e}odory class)}\\
\end{split}
\]

\section{Toeplitz matrices and de Branges spaces}
\label{sec:toeplitz}
We shall say that a $p\times 2p$ entire mvf $\gE=\begin{bmatrix}E_-&E_+\end{bmatrix}$ is a de Branges matrix with respect to $\DD$ if 
\begin{equation}
\label{eq:dec12a23}
\det E_+(\lambda)\ne 0\quad\textrm{for}\ \lambda\in{\DD} \quad\textrm{and}\quad E_+^{-1}E_-\in{\cS}_{\rm in}^{\ptp}.
\end{equation}
With each such mvf $\gE$, the space 
\begin{equation}
\label{eq:dec12b23}
{\cB}(\gE) =\{\textrm{entire $p\times 1$ vvf's}\ f:\ E_+^{-1}f\in H_2^p\ominus E_+^{-1}E_- H_2^p\}
\end{equation}
is a RKHS with RK  
\begin{equation}
		\label{eq:dec15i21}
		K_\omega(\lambda)=\frac{E_+(\lambda)E_+(\omega)^*-E_-(\lambda)E_-(\omega)^*}{1-\lambda\ol{\omega}}\quad 
		\textrm{for  $\lambda,\omega\in\CC$,  $\lambda\ol{\omega}\ne 1$}
	\end{equation}
and inner product
\begin{equation}
\label{eq:dec12c23}
\langle f,g\rangle_{{\cB}(\gE)}=\frac{1}{2\pi}\int_0^{2\pi}g(e^{it})\Delta_{\gE}(t)f(e^{it}) dt\quad \textrm{with}\quad \Delta_{\gE}(t)=(E_+(e^{it})E_+(e^{it})^*)^{-1}.
\end{equation}
Moreover, 
\begin{equation}
\label{eq:mar21a24}
\textrm{if}\quad {\gF}=\begin{bmatrix}{F_-}&{F_+}\end{bmatrix}=\begin{bmatrix}E_-&E_+\end{bmatrix}M\quad\textrm{and}\quad Mj_pM^*=j_p, \quad\textrm{then}\quad {\cB}(\gE)={\cB}({\gF}).
\end{equation}
In this section we shall show that the RKHS ${\cH}_G$ based on a positive definite matrix $G$ is a de Branges space with respect to $\DD$ if and only if $G$ is a block Toeplitz matrix. The proof rests upon the following characterization of de Branges spaces with respect to $\DD$:
\begin{thm}
	\label{thm:dec15b21}
	Let $\cH$ be a RKHS of entire $p\times 1$ vvf's 
	with RK $K_\omega(\lambda)$ on $\CC\times \CC$. Suppose further that there exists a nonzero point $\alpha\in\DD$  with  reflection $\alpha_\circ=1/\ol{\alpha}$ about $\TT$ such that 
	\begin{equation}
		\label{eq:dec15h21}
		K_\alpha(\alpha)\succ 0\quad and \quad K_{\alpha_\circ}(\alpha_\circ)\succ 0. 
	\end{equation}
	Then  ${\cH}$ is a de Branges space ${\cB}(\gE)$ with respect to $\DD$ 
	if and only if 
	\begin{enumerate}
		\item[\rm(1)] $R_\alpha{\cH}_\alpha\subseteq\cH$ and $R_{\alpha_\circ}{\cH}_{\alpha_\circ}\subseteq\cH$.
		\item[\rm(2)] The linear transformation 
		\begin{equation}
			\label{eq:dec15j21}
			S_\alpha=-\ol{\alpha}I+(1-\vert\alpha\vert^2)R_\alpha.
		\end{equation}
		maps $\cH_\alpha$ isometrically onto $\cH_{\alpha_\circ}$. 
	\end{enumerate}
	Moreover, in this case the blocks 
	$E_\pm(\lambda)$ of $\gE(\lambda)$ may be specified by the formulas
	\begin{equation}
		\label{eq:dec15k21}
		E_+(\lambda)=\rho_\alpha(\lambda)K_\alpha(\lambda)\{\rho_\alpha(\alpha)K_\alpha(\alpha)\}^{-1/2}	\end{equation}
	and
	\begin{equation}
		\label{eq:dec15m21}
		E_-(\lambda)=-\rho_{\alpha_\circ}(\lambda)K_{\alpha_\circ}(\lambda)\{-\rho_{\alpha_\circ}(\alpha_\circ)K_{\alpha_\circ}(\alpha_\circ)\}^{-1/2},
	\end{equation}
	where
	\begin{equation}
	\label{eq:nov17g23}
	\rho_\omega(\lambda)=1-\lambda\ol{\omega}.
\end{equation}
	This choice of 
	$\gE(\lambda)=\begin{bmatrix}E_-(\lambda)&E_+(\lambda)\end{bmatrix}$ is unique up to a multiplicative constant matrix on the right 
	that is $j_p$-unitary.
\end{thm} 

\begin{proof} See Theorem 5.2 in \cite{d23a} (which is formulated for a more general setting). \end{proof}

\begin{thm}
\label{thm:feb25b24}
Fix $\alpha\in\DD$ and let $G\in\CC^{\mtm}$ be any positive definite matrix. Then the isometry 
\begin{equation}
\label{eq:feb25c24}
\langle S_\alpha f,S_\alpha f\rangle_{{\cH}_G}=\langle f, f\rangle_{{\cH}_G}\quad\text{holds for every}\ f\in({\cH}_G)_\alpha
\end{equation}
if and only if $G$ is block Toeplitz.
\end{thm}

\begin{proof}
 If 
	$f=F u$ belongs to $({\cH}_G)_\alpha$ for some point $\alpha\in\DD$,  then 
$$
S_\alpha f=F(A-\ol{\alpha}I_m)v \quad\textrm{with} \quad v=(I_m-\alpha A)^{-1}u,\quad f=F(I_m-\alpha A)v \quad\textrm{and}\quad 	{\gee}_0^*v=F(\alpha)u=0.
$$
The  condition  $\langle S_\alpha f,S_\alpha f\rangle_{{\cH}_G} =\langle f,f\rangle_{{\cH}_G}$ holds for such $f$ if and only if

\[
\begin{split}
&\langle F(A-\ol{\alpha}I_m)v,  F(A-\ol{\alpha}I_m)v\rangle_{{\cH}_G}=
\langle F(I_m-\alpha A)v,  F(I_m-\alpha A)v\rangle_{{\cH}_G}, 	\\	
&\iff 
\langle G(A-\ol{\alpha}I_m)v,  (A-\ol{\alpha}I_m)v\rangle_{\rm st}=
\langle G(I_m-\alpha A)v,  (I_m-\alpha A)v\rangle_{\rm st}	\\	
&\iff v^*(G-A^*GA)v=0\quad\textrm{for all vectors $v\in\CC^m$ with ${\gee}_0^*v=0$.}
\end{split}
\] 
 But, this last condition holds if and only if 
	\begin{align*}
		[G-A^\ast GA]_{[1,n]}=0\iff g_{ij}=g_{i+1, j-1} \text{ for } i, j=0,\cdots, n 
		\end{align*}
i.e., if and only if $G$ is a block Toeplitz matrix.	
\end{proof}

\begin{thm}\label{thm:nov27c23}
	$\mathcal{H}_G$ is a de Branges space  ${\cB}(\gE)$ with respect to $\mathbb{D}$ if and only if $G$ is a block Toeplitz matrix. Moreover, in this case the matrix polynomials $E_\pm$ may be specified by the formulas 
	\begin{equation}
	\label{eq:nov27g23}
		  E_{+} (\lambda)= \frac{1-\lambda\ol{\alpha}}{(1-\vert\alpha\vert^2)^{1/2}}F(\lambda){\gz}_\alpha 
\quad \text{ and } \quad 
		 E_{-} (\lambda)=\frac{\lambda-\alpha}{(1-\vert\alpha\vert^2)^{1/2}}F(\lambda){\gz}_{\alpha_\circ}
	\end{equation}
	where $\alpha_\circ=1/\ol{\alpha}$ is the reflection of $\alpha$ with respect to $\TT$ and $\alpha$ is any point 
	in $\DD\setminus\{0\}$. 
\end{thm}
\begin{proof}
	The proof is an immediate corollary of  Theorems \ref{thm:dec15b21} and \ref{thm:feb25b24}. 
\end{proof}

	The next theorem presents more attractive formulas for $E_\pm$ that only involve two block columns of $\Gamma$. 
	
\begin{thm}
\label{thm:nov27a23}	
If the positive definite matrix $G$ is block Toeplitz, then the RK for ${\cH}_G$ 
 \begin{align}
 \label{eq:nov27i23}
		K_\omega(\lambda)=F(\lambda)\Gamma F(\omega)^\ast=\frac{E_+(\lambda)E_+(\omega)^*- E_-(\lambda)E_-(\omega)^*}{\rho_\omega(\lambda)} \quad (\lambda \overline{\omega}\neq 1),
	\end{align}
can be expressed 	in terms of the matrix polynomials 
	\begin{align}\label{eq:dec12d23}
		& E_+ (\lambda)= F(\lambda) {\gu}_0\quad  \text{and} \quad 
		E_ - (\lambda)=\wh{F}(\lambda) {\gu}_n.
	\end{align}
Moreover, 
	\begin{align}\label{eq:nov27c23}
				\frac{1}{2\pi} \int_0^{2\pi} v^\ast F(e^{i\theta})^\ast \Delta_{\gE}(e^{i\theta}) F(e^{i\theta})u \, d\theta =v^\ast G u,
	\end{align}
	for $\Delta_{\gE}=(E_+E_+^*)^{-1}$ and every choice of $u, v\in \mathbb{C}^{m}$.
\end{thm}	
 
\begin{proof}	Let $E_\pm^{(\alpha)}$ denote the mvf's defined in \eqref{eq:nov27g23}. 
Theorem \ref{thm:nov27c23} guarantees that formula \eqref{eq:dec15i21} holds with  $E_\pm=E_\pm^{(\alpha)}$. 
Formula \eqref{eq:nov27i23} is obtained by noting that the left hand side of \eqref{eq:dec15i21}  is independent of $\alpha$, whereas
$$
\lim_{\alpha\rightarrow 0}E_+^{(\alpha)}(\lambda) E_+^{(\alpha)}(\omega)^*=F(\lambda){\gu}_0 {\gu}_0^*F(\omega)^*
$$
and
$$
\lim_{\alpha\rightarrow 0}E_-^{(\alpha)}(\lambda) E_-^{(\alpha)}(\omega)^*=\wh{F}(\lambda){\gu}_n {\gu}_n^*\wh{F}(\omega)^*.
$$

Finally, to verify \eqref{eq:nov27c23}, let $v$ be any vector in $\CC^m$ and let $\omega_0,\ldots,\omega_n$ be any set of $n+1$ distinct points in $\CC$. Then,
since  the block companion matrix $\begin{bmatrix}F(\omega_0)^*&\cdots&F(\omega_n)^*\end{bmatrix}$ is invertible,
there exists 
 a set of vectors $x_0,\ldots,x_n\in\CC^p$ such that  $v=\Gamma\sum_{j=0}^n F(\omega_j)^*x_j$. Therefore,
 $Fv=\sum_{j=0}^n K_{\omega_j }x_j$ and
$$
\langle Fu,Fv\rangle_{{\cB}({\gE})}=\sum_{j=0}^n \langle Fu,K_{\omega_j} x_j\rangle_{{\cB}({\gE})}=\sum_{j=0}^n x_j^*F(\omega_j)u=v^*Gu
$$
for every vector $u\in\CC^m$.  In view of \eqref{eq:dec12c23}, the formula \eqref{eq:nov27c23} follows. 
\end{proof}

\begin{cor}
\label{cor:dec12a23}
If $G$ is a positive definite Toeplitz matrix with entries $g_{ij}=h_{i-j}\in\CC^{\ptp}$ for $i,j=0,\ldots,n$ and if $E_+$ is specified by formula \eqref{eq:dec12d23} and $\Delta_{\gE}(t)=(E_+(e^{it})E_+(e^{it})^*)^{-1}$, then \begin{equation}
\label{eq:oct24c23}
\frac{1}{2\pi}\int_0^{2\pi} e^{-ikt}\Delta_{\gE}(t)dt=h_k\quad\textrm{for}\ k=0,\pm 1,\ldots,\pm n.
\end{equation}
Thus, the density $\sigma(t)=\int_0^t\Delta_{\gE}(s)ds$ is a solution 
of the truncated trigonometric moment problem for $h_{-n},\ldots,h_n$.
\end{cor}
\begin{proof}
This is immediate from formula \eqref{eq:nov27c23} 
by choosing $u= {\gee}_i$ and $v= {\gee}_j$ with $i,j=0,\ldots,n$ and $i-j=k$.
\end{proof}

\begin{rem}
\label{rem:nov27b23}
The proof of Theorem \ref{thm:nov27c23} can be adapted to show that operator $R_0$ maps the RKHS $(\mathcal{H} _G)_0$ isometrically onto the RKHS 
$$
(\mathcal{H}_G)_\bullet=\{ R_0f: f\in \mathcal{H}_0\} =\left \{F\begin{bmatrix}
				v\\
				0
			\end{bmatrix} : v \in \mathbb{C}^{np} \right\}, 
$$
  if and only if $G$ is a block Toeplitz matrix. Moreover, the RK's of these two spaces
can be expressed as
$$		
			K_\omega^{(0)} (\lambda)
			=F(\lambda)(\Gamma-{\gu}_0{\gu}_0^*)F(\omega)^*=F(\lambda)\Gamma F(\omega)^*-E_+(\lambda)E_+(\omega)^*
$$
and		
$$		
			K_\omega^{(\bullet)} (\lambda)=   
			F(\lambda)(\Gamma-{\gu}_n{\gu}_n^*)F(\omega)^*=F(\lambda)\Gamma F(\omega)^*-E_-(\lambda)E_-(\omega)^*,
$$
respectively, with $E_\pm$ as in \eqref{eq:dec12d23}.
\end{rem}

\begin{rem}
\label{rem:dec12a23}
The formulas in \eqref{eq:dec12d23} were obtained by other means in \cite{d88} under the less restrictive assumptions 
that $G$ is an invertible Hermitian matrix such that $\gamma_{nn}^{(n)}$ is positive definite.
\end{rem}

\section{Hankel matrices and de Branges spaces}
\label{sec:hankel}
We shall say that a $p\times 2p$ entire mvf $\gE=\begin{bmatrix}E_-&E_+\end{bmatrix}$ is a de Branges matrix with respect to $\CC_+$ if 
\begin{equation}
\label{eq:dec12e23}
\det E_+(\lambda)\ne 0\quad\textrm{for}\ \lambda\in{\CC_+} \quad\textrm{and}\quad E_+^{-1}E_-\in{\cS}_{\rm in}^{\ptp}(\CC_+).
\end{equation}
With each such mvf $\gE$, the space 
\begin{equation}
\label{eq:dec12f23}
{\cB}(\gE) =\{\textrm{entire $p\times 1$ vvf's}\ f:\ E_+^{-1}f\in H_2^p\ominus E_+^{-1}E_- H_2^p\}
\end{equation}
is a RKHS with RK  
\begin{equation}
		\label{eq:dec15ii21}
		K_\omega(\lambda)=\frac{E_+(\lambda)E_+(\omega)^*-E_-(\lambda)E_-(\omega)^*}{-2\pi i(\lambda-\ol{\omega})}\quad 
		\textrm{for  $\lambda,\omega\in\CC$,  $\lambda\ne \ol{\omega}$}
	\end{equation}
and inner product
\begin{equation}
\label{eq:dec12g23}
\langle f,g\rangle_{{\cB}(\gE)}=\int_{-\infty}^\infty g(\mu)^*\Delta_{\gE}(\mu)f(\mu) d\mu\quad \textrm{with}\quad \Delta_{\gE}(t)=(E_+(\mu)E_+(\mu)^*)^{-1}.
\end{equation}
Moreover, the implication \eqref{eq:mar21a24} holds in this setting also.

In this section we shall show that the RKHS ${\cH}_G$ based on a positive definite matrix $G$ is a de Branges space with respect to $\CC_+$ if and only if $G$ is a block Hankel matrix. The proof rests upon the following characterization of de Branges spaces with respect to $\CC_+$:
\begin{thm}
	\label{thm:sep23b20}
	Let $\cH$ be a RKHS of entire $p\times 1$ vvf's 
	with RK $K_\omega(\lambda)$ on $\CC\times\CC$ such that there exists a point $\alpha\in\CC_+$  for which  
	\begin{equation}
		\label{eq:oct14d20}
		K_\alpha(\alpha)\succ 0\quad and \quad K_{\ol{\alpha}}(\ol{\alpha})\succ 0. 
	\end{equation}
	Then  there exists a $p\times 2p$ entire mvf $\gE(\lambda)=\begin{bmatrix}E_-(\lambda)&E_+(\lambda)\end{bmatrix}$       such that 
	\begin{equation}
		\label{eq:sep27a20}
		K_\omega(\lambda)=\frac{E_+(\lambda)E_+(\omega)^*-E_-(\lambda)E_-(\omega)^*}{-2\pi i(\lambda-\ol{\omega})}\quad \text{for\ $\lambda,\omega\in\CC$,\  $\lambda\ne\ol{\omega}$}
	\end{equation}
	if and only if 
	\begin{enumerate}
		\item[\rm(1)] $R_\alpha{\cH}_\alpha\subseteq\cH$ and $R_{\ol{\alpha}}{\cH}_{\ol{\alpha}}\subseteq\cH$.
		\item[\rm(2)] The linear transformation 
		\begin{equation}
			\label{eq:sep23c20}
			T_\alpha=I+(\alpha-\ol{\alpha})R_\alpha 
		\end{equation}
		maps $\cH_\alpha$ isometrically onto $\cH_{\ol{\alpha}}$. 
	\end{enumerate}
	Moreover, in this case the blocks 
	$E_\pm(\lambda)$ of $\gE(\lambda)$ may be specified by the formulas 	\begin{equation}
		\label{eq:sep4d20}
		E_+(\lambda)=\rho_\alpha(\lambda)K_\alpha(\lambda)\{\rho_\alpha(\alpha)K_\alpha(\alpha)\}^{-1/2}
	\end{equation}
	and
	\begin{equation}
		\label{eq:sep4e20}
		E_-(\lambda)=\rho_{\ol{\alpha}}(\lambda)K_{\ol{\alpha}}(\lambda)\{-\rho_{\ol{\alpha}}(\ol{\alpha})K_{\ol{\alpha}}(\ol{\alpha})\}^{-1/2}
	\end{equation}
	where
	\begin{equation}
	\label{eq:nov17ff23}
	\rho_\omega(\lambda)=-2\pi i(\lambda-\ol{\omega}).
	\end{equation}
 This choice of 
	$\gE(\lambda)=\begin{bmatrix}E_+(\lambda)&E_-(\lambda)\end{bmatrix}$ is unique up to a multiplicative  constant matrix on the right that is $j_p$-unitary.
	\end{thm} 

\begin{proof} See Theorem 7.1 in [DS]. \end{proof}

\begin{thm}
\label{thm:feb25d24}
Fix $\alpha\in\CC_+$ and let $G\in\CC^{\mtm}$ be any positive definite matrix. Then the isometry 
\begin{equation}
\label{eq:feb25c24}
\langle T_\alpha f,T_\alpha f\rangle_{{\cH}_G}=\langle f, f\rangle_{{\cH}_G}\quad\text{holds for every}\ f\in({\cH}_G)_\alpha
\end{equation}
if and only if $G$ is block Hankel.
\end{thm}

\begin{proof}
If $f=F u$ belongs to $({\cH}_G)_\alpha$ for some point $\alpha\in\CC_+$,  then 
$$
T_\alpha f=F(I_m-\ol{\alpha}A)v \quad\textrm{with} \quad v=(I_m-\alpha A)^{-1}u,\quad f=F(I_m-\alpha A)v \quad\textrm{and}\quad 	{\gee}_0^*v=F(\alpha)u=0.
$$
Thus, the  condition  $\langle T_\alpha f,T_\alpha f\rangle =\langle f,f\rangle$ holds for such $f$ if and only if

\[
\begin{split}
&\langle F(I_m-\ol{\alpha}A)v,  F(I_m-\ol{\alpha}A)v\rangle_{{\cH}_G}=
\langle F(I_m-\alpha A)v,  F(I_m-\alpha A)v\rangle_{{\cH}_G}, 	\\	
&\iff 
\langle G(I_m-\ol{\alpha}A)v,  (I_m-\ol{\alpha}A)v\rangle_{\rm st}=
\langle G(I_m-\alpha A)v,  (I_m-\alpha A)v\rangle_{\rm st}	\\	
&\iff v^*(GA-A^*G)v=0\quad\textrm{for all vectors $v\in\CC^m$ with ${\gee}_0^*v=0$.}
\end{split}
\] 
 But, this last condition holds if and only if 
	\begin{align*}
		[GA-A^\ast G]_{[1,n]}=0,
	\end{align*}
	which is, equivalent to $g_{ij} = g_{i+1,j-1}$ for $i = 0,\cdots, n$ and $j=1,\cdots, n$. Thus, the isometry holds if and only if G is a block Hankel
	matrix.
\end{proof}
\begin{thm}\label{thm: main Hankel}
The RKHS	$\mathcal{H}_G$ is a de Branges space  ${\cB}(\gE)$ with respect to $\bC_+$ if and only if $G$ is a block Hankel matrix. Moreover, in this case the blocks $E_\pm$ of $\gE$ may be specified by the formulas 
\begin{equation}
\label{eq:oct29a23}	
E_+(\lambda)=\rho_\alpha(\lambda)F(\lambda)\Gamma F(\alpha)^*	(\rho_\alpha(\alpha)F(\alpha)\Gamma F(\alpha)^*)^{-1/2}=\frac{\rho_\alpha(\lambda)}{\sqrt{\rho_\alpha(\alpha)}}F(\lambda){\gz}_{\alpha}
\end{equation}
and
\begin{equation}
\label{eq:oct29b23}	
E_-(\lambda)=\rho_{\ol{\alpha}}(\lambda)F(\lambda)\Gamma F(\ol{\alpha})^*	(\rho_\alpha(\alpha)F(\ol{\alpha})\Gamma F(\ol{\alpha})^*)^{-1/2}=\frac{\rho_{\ol{\alpha}}(\lambda)}{\sqrt{\rho_\alpha(\alpha)}}F(\lambda){\gz}_{\ol{\alpha}}
\end{equation}	
for any choice of the point $\alpha\in\CC_+$. In this setting, $\rho_\omega(\lambda)=-2\pi i(\lambda-\ol{\omega})$.
\end{thm}
\begin{proof}
	The proof is an immediate consequence of Theorems \ref{thm:sep23b20} and \ref{thm:feb25d24}.
\end{proof}

\begin{cor}
\label{cor:dec12b23} 
If $G$ is a positive definite block Hankel matrix with entries $g_{ij}=h_{i+j}\in\CC^{\ptp}$ for $i,j=0,\ldots,n$ and if $E_+$ is specified by formula \eqref{eq:oct29a23} and $\Delta_{\gE}=(E_+E_+^*)^{-1}$, then 
\begin{equation}
\label{eq:oct24b23}
\int_{-\infty}^\infty \mu^k\Delta_{\gE}(\mu)d\mu=h_k\quad\textrm{for}\ k=0,\ldots,2n.
\end{equation}
Thus,  $\sigma(\mu)=\int_{-\infty}^\mu\Delta_{\gE}(x)dx$ is a solution 
of the truncated Hamburger matrix moment problem for $h_0,\ldots,h_{2n}$..
\end{cor}

\begin{proof}
This is immediate from  the formula (which follows from \eqref{eq:dec12g23}) 
\[
\int_{-\infty}^{\infty} v^\ast F(\mu)^\ast \Delta_{\gE}(\mu) F(\mu)u \, d\mu =v^\ast G u,
\]
by choosing $u= {\gee}_i$ and $v= {\gee}_j$ with $i,j=0,\ldots,n$ and $i+j=k$.
\end{proof}
 
\begin{rem}
\label{rem:dec12b23}
Formulas for $E_\pm$ that only involve two block columns of $\Gamma$ may be obtained by setting $\alpha=re^{i\theta}$ and letting $r\uparrow \infty$. Although such formulas are attractive for matrix computations, they do not simplify the 
analysis needed to obtain solutions to the truncated matrix Hamburger moment problem. 
Such  formulas were obtained by other means in Section 9 of \cite{d89a}
 under the less restrictive assumptions 
that $G$ is an invertible Hermitian matrix such that $\gamma_{nn}^{(n)}$ is positive definite. 
\end{rem}

\section{Characterizations of de Branges spaces and matrix identities}
\label{sec:dbspacecharacterizations}
This section is devoted to a brief survey of the ideas behind the characterizations of the two classes of de Branges spaces presented in Theorems \ref{thm:sep23b20} and \ref{thm:dec15b21}, and subsequently to exploit the two basic isometries $S_\alpha$ of \eqref{eq:dec15j21} and $T_\alpha$ of \eqref{eq:sep23c20} to obtain a number of matrix identities that include the well known formulas of Gohberg-Heinig and Gohberg-Semencul. 
 
 \subsection{The Hankel setting} 
 Theorem \ref{thm:sep23b20} is an extension to RKHS's of vector valued entire functions of a very beautiful observation of L. de Branges that was formulated for scalar valued entire functions in Sections 22 and 23 of 
\cite{br68a} that deserves to be better known. To begin with, the operator $T_\alpha$, $\alpha\in\CC_+$, which may look a little strange at first glance, is just multiplication  by the inverse of an elementary Blaschke product:
\begin{equation}
\label{eq:nov20a23}
 T_\alpha : f\in{\cH}_\alpha\mapsto b_\alpha^{-1}f, \quad\textrm{where }\ b_\alpha(\lambda)=\frac{\lambda-\alpha}{\lambda-\ol{\alpha}}.
\end{equation}
It is written in the form \eqref{eq:sep23c20}
 in order to express the operator in terms of $R_\alpha$. In particular, 
\begin{equation}
\label{eq:nov22f23}
\langle T_\alpha f,T_\alpha f\rangle_{\cH}-\langle f,f\rangle_{\cH} =
\langle R_\alpha f, f\rangle_{\cH} -\langle f, R_\alpha f \rangle_{\cH} -(\alpha-\overline{\alpha})\langle R_\alpha f, R_\alpha f \rangle_{\cH},
\end{equation}
which serves to relate the condition \eqref{eq:feb25c24} to another fundamental identity due to de Branges; see \eqref{eq:jun23a24} below. 
Moreover, since $g\in{\cH}_{\ol{\alpha}}\implies g=T_\alpha f$ for exactly one vvf $f\in{\cH}_\alpha$ and $T_{\ol{\alpha}}g=T_{\ol{\alpha}}T_\alpha f=f$, it is readily checked that $T_{\ol{\alpha}}$ maps ${\cH}_{\ol{\alpha}}$ isometrically onto ${\cH}_\alpha$.

In particular,
$$
v^*b_\alpha(\lambda)^{-1}K_\omega^{(\alpha)}(\lambda)u=\langle T_\alpha K_\omega^{(\alpha)}u, K_\lambda^{(\ol{\alpha})}v\rangle_{\cH}=\langle K_\omega^{(\alpha)}u, T_{\ol{\alpha}}K_\lambda^{(\ol{\alpha})}v\rangle_{\cH}=\{u^*(T_{\ol{\alpha}}K_\lambda^{(\ol{\alpha})}(\omega)v\}^*
$$
for every pair of points $\lambda,\omega\in\CC$ and every pair of vectors $u,v\in\CC^p$, which upon invoking the 
formula \eqref{eq:nov22e23} 
for $K_\omega^{(\alpha)}(\lambda)$ and its counterpart for $K_\omega^{(\ol{\alpha})}(\lambda)$ yields the identity
\begin{equation}
\label{eq:nov20b23}
\frac{\lambda-\ol{\alpha}}{\lambda-\alpha}\left\{K_\omega(\lambda)-K_\alpha(\lambda)K_\alpha(\alpha)^{-1}K_\omega(\alpha)\right\}=\frac{\ol{\omega}-\ol{\alpha}}{\ol{\omega}-\alpha}\left\{
 K_\omega(\lambda)-K_{\ol{\alpha}}(\lambda)K_{\ol{\alpha}}(\ol{\alpha})^{-1}K_\omega(\ol{\alpha})\right\}.
\end{equation}
Formulas \eqref{eq:sep4d20} and \eqref{eq:sep4e20} emerge upon solving for $K_\omega(\lambda)$ in \eqref{eq:nov20b23}, if you have enough courage and stamina to push through to the end. Moreover, there are other significant equivalent matrix identities:

\begin{thm}
\label{thm:feb27a24}
If $G$ is positive definite and $\alpha\in\CC\setminus\RR$, 
 then the following statements are equivalent. (Moreover, the  identities (3), (4), (5) are equivalent for every choice of $\alpha,\lambda,\omega\in\CC$.)
\begin{enumerate}
\item[\rm(1)] $G$ is a block Hankel matrix. 
\item[\rm(2)] The isometry $\langle T_\alpha f,T_\alpha f\rangle_{{\cH}_G}=\langle  f,f\rangle_{{\cH}_G}$ holds for every $f\in({\cH}_G)_\alpha$ and $\alpha\in\CC\setminus\RR$. 
\item[\rm(3)] $(\lambda-\ol{\alpha})(\ol{\omega}-\alpha) F(\lambda)N_\alpha F(\omega)^*=
(\lambda-\alpha)(\ol{\omega}-\ol{\alpha} )F(\lambda)N_{\ol{\alpha}} F(\omega)^*$.
\item[\rm(4)] $(\lambda-\ol{\alpha}) F(\lambda)N_\alpha (I_m-\alpha A^*)=
(\lambda-\alpha)F(\lambda)N_{\ol{\alpha}} (I_m-\ol{\alpha}A^*)$.
\item[\rm(5)] $(I_m-\ol{\alpha}A) N_\alpha (I_m-\alpha A^*)=
(I_m-\alpha A)N_{\ol{\alpha}} (I_m-\ol{\alpha}A^*)$.
\end{enumerate}
\end{thm}

\begin{proof} The equivalence of (1) and (2)  has already been established in Theorem 
\ref{thm:feb25d24}. The formula in (3) is equivalent to the identity  \eqref{eq:nov20b23}, since \eqref{eq:nov22e23} ensures that $K_\omega^{(\beta)}(\lambda)=F(\lambda)N_\beta F(\omega)^*$ for every point $\beta\in\CC$, and in particular for $\beta=\alpha$ and $\beta=\ol{\alpha}$. Thus, $(3)\implies(1)$, since 
  the proof that the RKHS with RK $K_\omega(\lambda)$ is a de Branges  space, rests on \eqref{eq:nov20b23}; see e.g., \cite{ds17} for details, though the basic  idea is due to de Branges \cite{br68b}.

Moreover, 
$(3)\implies (4)$, since
\[
\begin{split}
N_\alpha F(\omega)^*&=N_\alpha\{F(\omega)-F(\alpha)\}^*=(\ol{\omega}-\ol{\alpha})N_\alpha (R_\alpha F)(\omega)^*\\
&=(\ol{\omega}-\ol{\alpha})N_\alpha (I_m-\ol{\alpha}A^*)^{-1}A^*F(\omega)^*
\end{split}
\]
and, by the same reasoning, but with $\ol{\alpha}$ in place of $\alpha$, 
\[
N_{\ol{\alpha}} F(\omega)^*=(\ol{\omega}-\alpha)N_{\ol{\alpha}} (I_m-\alpha A^*)^{-1}A^*F(\omega)^*
\]
Thus, in view of (3),  
\[
(\lambda-\ol{\alpha})F(\lambda)N_\alpha (I_m-\ol{\alpha}A^*)^{-1}A^*F(\omega)^*=
(\lambda-\alpha)F(\lambda)N_{\ol{\alpha}} (I_m-{\alpha}A^*)^{-1}A^*F(\omega)^*
\] 
for every $\omega\in\CC$ and hence
\[
(\lambda-\ol{\alpha})F(\lambda)N_\alpha (I_m-\ol{\alpha} A^*)^{-1}A^*=
(\lambda-\alpha)F(\lambda)N_{\ol{\alpha}} (I_m-\alpha A^*)^{-1}A^*.
\] 
But, as $A^*A+{\gee}_0{\gee}_0^*=I_m$ and $N_\alpha (I_m-\ol{\alpha}A^*)^{-1}{\gee}_0=N_\alpha F(\alpha)^*=0$, the last identity identity implies (4), i.e, $(3)\implies (4)$. 

The proof that $(4)\implies (5)$ is similar to the proof that $(3)\implies (4)$. It rests on the observation that 
$$
F(\lambda)N_\alpha=\{F(\lambda)-F(\alpha)\}N_\alpha=(\lambda-\alpha) (R_\alpha F)(\lambda)=(\lambda-\alpha) F(\lambda)A(I_m-\alpha A)^{-1}N_\alpha
$$
and another application of the identity $A^*A+{\gee}_0{\gee}_0^*=I_m$. It remains only to verify that $(5)\implies (3)$. Towards this end,  observe that (5) implies that
$$
A(I_m-\alpha A)^{-1}N_\alpha (I_m-\ol{\alpha}A^*)^{-1}A^*=A(I_m-\ol{\alpha} A)^{-1}N_{\ol{\alpha}} (I_m-{\alpha}A^*)^{-1}A^*
$$
and hence that 
$$
F(\lambda)A(I_m-\alpha A)^{-1}N_\alpha (I_m-\ol{\alpha}A^*)^{-1}A^*F(\omega)^*=
F(\lambda)A(I_m-\ol{\alpha} A)^{-1}N_{\ol{\alpha}} (I_m-{\alpha}A^*)^{-1}A^*F(\omega)^*.
$$
The proof is completed by invoking the formulas 
$$
F(\lambda)A(I_m-\beta A)^{-1}N_\beta =\frac{F(\lambda)-F(\beta)}{\lambda-\beta}N_\beta =\frac{F(\lambda)}{\lambda-\beta}N_\beta
$$
and
$$
N_\beta (I_m-\ol{\beta}A^*)^{-1}A^*F(\omega)^*=N_\beta \left\{\frac{F(\omega)-F(\beta)}{\omega-\beta}\right\}^*
=N_\beta \left\{\frac{F(\omega)}{\omega-\beta}\right\}^*
$$
for $\beta=\alpha$ and $\beta=\ol{\alpha}$.
\end{proof}


\subsection{The Toeplitz setting} 
Theorem \ref{thm:dec15b21} is established by much the same strategy that was used to justify Theorem \ref{thm:sep23b20}. The key is to note that  the pair $\alpha\in\CC_+$ and its reflection $\ol{\alpha}$ about $\RR$ should be replaced by the pair $\alpha\in\DD\setminus\{0\}$ and its reflection $\alpha_\circ=1/\ol{\alpha}$ about $\TT$ and that $T_\alpha$ should be replaced by $S_\alpha$, alias the operator of multiplication by the inverse of the Blaschke product $b_\alpha(\lambda)=(\lambda-\alpha)/(1-\lambda \ol{\alpha})$ for the disk. In this setting, 
\begin{equation}
\label{eq:dec17a23}
\langle S_\alpha f,S_\alpha f\rangle_{\cH}-\langle f,f\rangle_{\cH}=
\langle f,f\rangle_{\cH}+\alpha\langle R_\alpha f, f\rangle_{\cH} +\ol{\alpha}\langle f, R_\alpha f \rangle_{\cH} -(1-\vert\alpha\vert^2)\langle R_\alpha f, R_\alpha f \rangle_{\cH}=0 
\end{equation}
for $f\in{\cH}_\alpha$. 
This is the restriction of an analogue of de Branges identity   that originates in the work of Ball and Rovnyak \cite{ba75}; see \eqref{eq:jun23b24} below.   For additional perspective on these identities and the connection
\begin{equation}
\label{eq:feb28a24}
\frac{1-\lambda\ol{\alpha}}{\lambda-\alpha}K_\omega^{(\alpha)}(\lambda)=\frac{\ol{\omega}-\ol{\alpha}}{1-\alpha\ol{\omega}}K_\omega^{(\alpha_\circ)}(\lambda)
\end{equation}
between the RK's $K_\omega^{(\alpha)}(\lambda)$ and $K_\omega^{(\alpha_\circ)}(\lambda)$ see  \cite{d23a}.

\begin{thm}
\label{thm:feb28a24}
If $G$ is positive definite and $\alpha\in\DD\setminus\{0\}$, 
 then the following statements are equivalent. (Moreover, the  identities (3), (4), (5) are equivalent for every choice of $\alpha\in\CC\setminus\{0\}$ and $\lambda,\omega\in\CC$.)
\begin{enumerate}
\item[\rm(1)] $G$ is a block Toeplitz matrix. 
\item[\rm(2)] The isometry $\langle S_\alpha f,S_\alpha f\rangle_{{\cH}_G}=\langle  f,f\rangle_{{\cH}_G}$ holds for every $f\in({\cH}_G)_\alpha$ and $\alpha\in\DD$. 
\item[\rm(3)] $(1-\lambda\ol{\alpha})(1-\alpha\ol{\omega})F(\lambda)N_\alpha F(\omega)^*=(\lambda-\alpha)(\ol{\omega}-\ol{\alpha})F(\lambda)N_{\alpha_\circ}F(\omega)^*$
\item[\rm(4)] $(1-\lambda\ol{\alpha})F(\lambda)N_\alpha (I_m-\ol{\alpha} A^*)^{-1}=-(\lambda-\alpha)F(\lambda)N_{\alpha_\circ}(\alpha I_m-A^*)^{-1}$.
\item[\rm(5)] $(\ol{\alpha}I_m-A)N_\alpha(\alpha I_m-A^*)=(I_m-\alpha A)N_{\alpha_\circ}(I_m-\ol{\alpha}A^*)$.
\end{enumerate}
\end{thm}

\begin{proof} 
The equivalence of (1) and (2) is established in Theorem \ref{thm:feb25b24}. The formula in (3) is clearly equivalent to \eqref{eq:feb28a24}, and hence also to (1) and (2), since 
  the proof that the RKHS with RK $K_\omega(\lambda)$ is a de Branges  space, rests on\eqref{eq:feb28a24}; see e.g., \cite{d23a} for details.

The remaining two identities are obtained  by multiple use of the observation that 
 that $F(\alpha)N_\alpha=0$ and hence that 
\begin{equation}
\label{eq:jun13a24}
F(\lambda)N_\alpha= (\lambda-\alpha)\frac{F(\lambda)-F(\alpha)}{\lambda-\alpha}N_\alpha=(\lambda-\alpha)F(\lambda)A(I_m-\alpha A)^{-1}N_\alpha.
\end{equation}

In particular, the implication $(3)\implies(4)$, then follows by substituting the formula
\[
\begin{split}
N_\alpha F(\omega)^*&=(\ol{\omega}-\ol{\alpha})N_\alpha \left\{
\frac{F(\omega)-F(\alpha)}{\omega-\alpha}\right\}^*
=(\ol{\omega}-\ol{\alpha})N_\alpha \{F(\omega)A(I_m-\alpha A)^{-1}\}^*\\
&=(\ol{\omega}-\ol{\alpha})N_\alpha (I_m-\ol{\alpha} A^*)^{-1}A^*F(\omega)^*
\end{split}
\]
and its analogue \[
N_{\alpha_\circ} F(\omega)^*=(\ol{\omega}-\ol{\alpha_\circ})N_{\alpha_\circ} (I_m-\ol{\alpha_\circ} A^*)^{-1}A^*F(\omega)^*
\]
for $\alpha_\circ=1/\ol{\alpha}$ into the formula in (3). This yields the formula
\[
(1-\lambda\ol{\alpha})F(\lambda)N_\alpha (I_m-\ol{\alpha}A^*)^{-1}A^*F(\omega)^*=-(1/\alpha)(\lambda-\alpha)(F(\lambda)N_{\alpha_\circ}(I_m-\ol{\alpha_\circ}A^*)^{-1}A^*F(\omega)^*.
\]
But, since this is valid for all points $\omega\in\CC$, the same formula is valid when the term $F(\omega)^*$ is removed and then both sides are multiplied on the right by $A$, i.e., 
\[
(1-\lambda\ol{\alpha})F(\lambda)N_\alpha (I_m-\ol{\alpha}A^*)^{-1}A^*A=-(\lambda-\alpha)F(\lambda)N_{\alpha_\circ}(\alpha I_m-A^*)^{-1}A^*A.
\]
The last identity reduces  to the formula in  (4),  because $F(\beta)N_\beta=0$ and hence 
\[
N_\beta (I_m-\ol{\beta}A^*)^{-1}A^*A=N_\beta (I_m-\ol{\beta}A^*)^{-1}({\gee}_0{\gee}_0^*+A^*A)=N_\beta (I_m-\ol{\beta}A^*)^{-1}
\]
for every point $\beta\in\CC$. 

The verification that $(4)\implies$ 
\begin{equation}
\label{eq:mar1a24}
A(I_m-\alpha A)^{-1}N_\alpha (I_m-\ol{\alpha} A^*)^{-1}A^*=A(\ol{\alpha}I_m-A)^{-1} N_{\alpha_\circ}(\alpha I_m-A^*)^{-1}A^* 
\end{equation}
is similar to the verification that $(3)\implies(4)$. It rests on the fact that $F(\alpha)N_\alpha=0$ and hence that 
 \eqref{eq:jun13a24} and its analogue for $\alpha_\circ$ hold. Then, since $I_m={\gee}_0{\gee}_0^*+A^*A$ and ${\gee}_0^*(I_m-\alpha A)^{-1}N_\alpha=0$, \eqref{eq:mar1a24} implies that 
$$
(I_m-\alpha A)^{-1}N_\alpha (I_m-\ol{\alpha} A^*)^{-1}=(\ol{\alpha}I_m-A)^{-1} N_{\alpha_\circ}(\alpha I_m-A^*)^{-1},
$$
which is equivalent to (5).

It remains to verify the implication $(5)\implies(3)$. Clearly, (5) implies \eqref{eq:mar1a24}, which implies that 
\[
\begin{split}
F(\lambda)A(I_m-\alpha A)^{-1}N_\alpha (I_m-\ol{\alpha} A^*)^{-1}A^*&=F(\lambda)A(\ol{\alpha}I_m-A)^{-1} N_{\alpha_\circ}(\alpha I_m-A^*)^{-1}A^*\\
&=\vert\alpha_\circ\vert^2 F(\lambda)A(I_m-\alpha_\circ A)^{-1} N_{\alpha_\circ}( I_m-\ol{\alpha_\circ}A^*)^{-1}A^*\\
\end{split}
\]
and hence, in view of \eqref{eq:nov17f23}, that
\[
\frac{F(\lambda)-F(\alpha)}{\lambda-\alpha}N_\alpha (I_m-\ol{\alpha} A^*)^{-1}A^*=\vert\alpha_\circ\vert^2 \frac{F(\lambda)-F(\alpha_\circ)}{\lambda-\alpha_\circ}N_{\alpha_\circ}( I_m-\ol{\alpha_\circ}A^*)^{-1}A^*.
\]
Thus, as $F(\alpha)N_\alpha=0$ and $F(\alpha_\circ)N_{\alpha_\circ}=0$, 
\[
(\lambda-\alpha_\circ)F(\lambda)N_\alpha (I_m-\ol{\alpha} A^*)^{-1}A^*=(\lambda-\alpha)\vert\alpha_\circ\vert^2 F(\lambda)N_{\alpha_\circ}( I_m-\ol{\alpha_\circ}A^*)^{-1}A^*.
\]
The formula in (3) is obtained by multiplying both sides on the right by $F(\omega)^*$ and invoking formula \eqref{eq:nov17f23} again.
\end{proof}

Since ${\gz}_0={\gu}_0$, $N_0=\Gamma-{\gu}_0{\gu}_0^*$. The supplementary {\bf notation}  $N_\bullet=\Gamma-{\gu}_n{\gu}_n^*$ plays the role of $N_{\alpha_\circ}$ when $\alpha=0$.

\begin{cor}
\label{cor:mar20a24}
If $G$ is block Toeplitz, then the following four equivalent identities are in force:
\begin{enumerate}
\item[\rm(1)] $F(\lambda)N_0 F(\omega)^*=\lambda\ol{\omega}F(\lambda)N_\bullet F(\omega)^*$.
\item[\rm(2)] $F(\lambda)N_0 A^*=\lambda F(\lambda)N_\bullet$.
\item[\rm(3)] $AN_0 A^\ast=N_\bullet$.
\item[\rm(4)] $AN_0=N_\bullet A$.
\end{enumerate}
\end{cor}

\begin{proof}
The first three identities follow from formulas $(3)-(5)$ in Theorem \ref{thm:feb28a24} by letting $\alpha\rightarrow 0$, since 
$$
\lim_{\alpha\rightarrow 0} N_\alpha=\Gamma-{\gu}_0{\gu}_0^*=N_0\quad\textrm{and}\quad  
\lim_{\alpha\rightarrow 0}N_{\alpha_\circ}=\Gamma-{\gu}_n{\gu}_n^*=N_\bullet.
$$
The equivalence of (3) and (4) rests on the fact that ${\gee}_0{\gee}_0^*+A^*A=I_m=AA^*+{\gee}_n{\gee}_n^*$ and
that ${\gee}_0^*N_0=0$ and ${\gee}_n^*N_\bullet=0$. Thus, if (3) is in force, then 
\[
N_\bullet A=AN_0A^*A=AN_0(A^*A+{\gee}_0{\gee}_0^*)=AN_0,
\]
i.e., (4) holds. Conversely, if (4) is in force, then
\[
AN_0A^*=N_\bullet AA^*=N_\bullet (AA^*+{\gee}_n{\gee}_n^*)=N_\bullet,
\]
i.e., (3)  holds.
\end{proof}

\begin{rem}
\label{rem:oct30a30}
The third formula  in Corollary \ref{cor:mar20a24}  is equivalent to the well known Gohberg-Heinig formula for the inverse of a Hermitian block Toeplitz matrix (see \cite{gh74})
\begin{equation}
	\label{eq:oct26a23}
	\Gamma=\sum_{j=0}^nA^j({\gu}_n{\gu}_n^*-A{\gu}_0{\gu}_0^*A^*)(A^*)^j.
	\end{equation}
		\end{rem}


\section{Solutions to the truncated matrix trigonometric moment problem}
\label{sec:solutionstrig}
In this section, 
 following  the strategy of \cite{mm03} and \cite{mm05}, we use the Arov-Grossman parameterization of minimal unitary extensions of isometric operators applied to  the isometric operator $R_0$ (alias $S_\alpha$ with $\alpha=0$), to obtain a description of the set 
\[
\Sigma_{\rm T}\quad \textrm{of solutions $\sigma$ to the truncated trigonometric matrix moment problem}.
\]
However, we use the underlying de Branges space structure to compute  projections in terms of appropriately chosen reproducing kernels rather than orthogonal bases. This serves to ease the bookkeeping and yields pleasing expressions for the principal formulas. This description is presented in Theorem \ref{thm:may7a24}. It is the main result of this section. Analogous calculations can be carried out for $S_\alpha$ with $\alpha\in\DD\setminus\{0\}$. The details are left to the reader.

Recall that when $G$ is block Toeplitz, then ${\cH}_G$ is a de Branges space ${\cB}({\gE})$ with ${\gE}=\begin{bmatrix}E_-&E_+\end{bmatrix}$,
\begin{equation}
\label{eq:mar10a24}
E_+ (\lambda)= 
F(\lambda) {\gu}_0\quad 
 \textrm{ and } \quad
		E_ - (\lambda)=\wh{F}(\lambda) {\gu}_n=\lambda F(\lambda) {\gu}_n
\end{equation}
and inner product \eqref{eq:dec12c23}, and  let
\begin{equation} 
\label{eq:mar10b24}
\Phi_{\gE}(\omega)=\frac{1}{2\pi}\int_0^{2\pi} \frac{e^{it}+\omega}{e^{it}-\omega}\Delta_{\gE}(t)dt,\quad\textrm{where}\quad 
\Delta_{\gE}(t)=(E_+(e^{it})E_+(e^{it})^*)^{-1}.
\end{equation}
The description of $\Sigma_T$ will be formulated in terms of a linear fractional transformation with coefficients that are expressed in terms of the entries $\Theta_{ij}$ of the $2p\times 2p$ matrix polynomial
\begin{equation}
\label{eq:mar10d24}
\Theta=\begin{bmatrix}\Theta_{11}&\Theta_{12}\\\Theta_{21}&\Theta_{22}\end{bmatrix}=
\frac{1}{\sqrt{2}}\begin{bmatrix}E_-^\circ&E_+^\circ\\ E_-&E_+\end{bmatrix},
\end{equation}
where $E_\pm^\circ$ 
are matrix polynomials that can be expressed in terms of the matrices

\begin{align}
	\mathbb{L}=\begin{bmatrix}
		& h_0 & 0 & \cdots & 0\\
		& 2h_1 & h_0& \cdots &0 \\
		&\vdots & & & \vdots& \\
		&2h_n& 2h_{n-1} & \cdots &h_0
	\end{bmatrix}\quad \text{and}\quad \LL^*=2G-\LL\end{align}
by the formulas
\begin{equation}
 \label{eq:nov30b23}
 E_+^\circ=F\LL{\gu}_0  \quad and \quad E_-^\circ=-\wh{F}\LL^*{\gu}_n,
 \end{equation}
respectively.


For ease of future reference we list the assumptions and definitions that will be in force throughout this section and in future discussion of the matrix trigonometric moment problem.
\begin{enumerate}
\item[\rm(T1)] $G$ is a positive definite block Toeplitz matrix and $\Gamma=G^{-1}$.
\item[\rm(T2)]  $\rho_\omega(\lambda)=1-\lambda\ol{\omega}$ and $f^\#(\lambda)=f(1/\ol{\lambda})^*$ for $\lambda\in\CC\setminus\{0\}$.
\item[\rm(T3)] The mvf's $E_\pm$, $\Phi_{\gE}$, $\Theta$ and $E_\pm^\circ$   are specified by formulas \eqref{eq:mar10a24}, 
 \eqref{eq:mar10b24}, \eqref{eq:mar10d24} and \eqref{eq:nov30b23}, 
 respectively, and the inner product in ${\cB}(\gE)$ is given by \eqref{eq:dec12c23} (and \eqref{eq:aug24c23}).
\end{enumerate}

Our first order of business is to find alternate expressions for the polynomials $E_\pm^\circ$ in terms of $E_\pm$ and $\Phi_{\gE}$ by exploiting the two different ways of evaluating 
$\frac{1}{2\pi}\int_0^{2\pi}\Delta_{\gE}(t)(R_\omega E_+)(e^{it})dt$ 
and  formula \eqref{eq:nov21c23} for computing the orthogonal projection of the columns of $\Delta_{\gE}E_\pm$ onto 
$H_2^p(\DD)$. We need a preliminary lemma.

\begin{lem}
\label{lem:may12a24}
If $(T1)$--$(T3)$ are in force, then:
\begin{enumerate}
\item[\rm(1)]  $\det\,E_+(\omega)\ne 0$ for $\vert\omega\vert\le 1$ and $\det\,\omega^{n+1}E_-^\#(\omega)\ne 0$
for $\vert\omega\vert\le 1$. 
\item[\rm(2)]  The mvf's $E_+(\omega)^{-1}$ and 
 $(\omega^{n+1}E_-^\#(\omega))^{-1}$ belong to the Wiener algebra $\mathcal{W}_{+}^{p\times p}(\mathbb{T})$ of $\ptp$ mvf's with summable Fourier coefficients. 
  \end{enumerate}
 \end{lem}
 
 \begin{proof}
 The asserted invertibility of $E_+$ follows from the fact that 
 \[
 E_+(\omega)E_+(\omega)^*-E_-(\omega)E_-(\omega)^*=(1-\vert\omega\vert^2)K_\omega(\omega)
 \]
since the RK  $K_\omega(\omega)=F(\omega)\Gamma F(\omega)^*$ is positive definite for every point $\omega\in\CC$.
 Thus, the reverse polynomial
$$
	E_+^\flat(\lambda)={\gee}_0^*(I_m-\lambda A)^{-1}(Z\Gamma Z){\gee}_0\{{\gee}_0^*Z\Gamma Z{\gee}_0\}^{-1/2}=
	{\gee}_0^*(I_m-\lambda A)^{-1} (Z\Gamma) {\gee}_n\{{\gee}_n^*\Gamma {\gee}_n\}^{-1/2}
	$$ 
	that is based on the positive definite block Toeplitz matrix $ZGZ$ with $Z=\begin{bmatrix}{\gee}_n&\cdots&{\gee}_0\end{bmatrix}$ is also invertible when $\vert\lambda\vert\le 1$. Consequently, the identity 
\[
 \begin{split}
 E_+^\flat(\lambda)&=
   \sum_{j=0}^n \lambda^j \gamma_{n-j,n}^{(n)} (\gamma_{nn}^{(n)})^{-1/2}=
  \sum_{k=0}^n \lambda^{n-k} \gamma_{kn}^{(n)} (\gamma_{nn}^{(n)})^{-1/2}\\ 
  &=\lambda^n F(1/\lambda){\gu}_n=\lambda^{n+1}\wh{F}(1/\lambda){\gu}_n=\lambda ^{n+1}E_-(1/\lambda), 
    \end{split}
 \]  
(which 
is equivalent to the formula $E_{-}^\# (\lambda)^{-1}=\lambda ^{n+1} (E_ +^\flat(\ol{\lambda})^*)^{-1}$) 
ensures that $\lambda^{n+1}E_-^\#(\lambda)$ is invertible for $\vert\lambda\vert\le 1$. 
This completes the proof of (1); (2) then follows from a well known theorem of Wiener.
\end{proof}

\begin{thm}
\label{thm:may8a24}
If $T(1)$--$(T3)$ are in force, then
\begin{equation}
\label{eq:may9d24}
\begin{split}
 E_+^\circ(\omega)&=h_0E_+(\omega)-2\omega {\gee}_0^*GA(I_m-\omega A)^{-1}{\gu}_0
  =h_0E_+(\omega)-\frac{2\omega}{2\pi}\int_0^{2\pi}\Delta_{\gE}(t)(R_\omega E_+)(t)dt\\ &
  =\Phi_{\gE}(\omega)E_+(\omega)
 \end{split}
 \end{equation}
and
\begin{equation}
\label{eq:may7a24}
\begin{split}
E_-^\circ(\omega)&=h_0E_-(\omega)-2\omega{\gee}_0^*G(I_m-\omega A)^{-1}{\gu}_n
=h_0E_-(\omega)-\frac{2\omega}{2\pi}\int_0^{2\pi}\Delta_{\gE}(t)(R_\omega E_-)(e^{it})dt\\ 
&=\Phi_{\gE}(\omega)E_-(\omega) -2E_-^\#(\omega)^{-1}.
\end{split}
\end{equation}
for $\omega\in\DD$.
\end{thm}

\begin{proof} The proof is broken into a number of steps.
\bigskip

\noindent
{\bf 1.}\ {\it Verification of the formulas}
\begin{equation}
\label{eq:may10a24}
\frac{1}{2\pi}\int_0^{2\pi}\Delta_{\gE}(t) (R_\omega F)(e^{it})dt={\gee}_0^*GA(I_m-\omega A)^{-1}
\end{equation}
and
\begin{equation}
\label{eq:may10b24}
\frac{1}{2\pi}\int_0^{2\pi}\Delta_{\gE}(t) (R_\omega \wh{F})(e^{it})dt={\gee}_0^*G(I_m-\omega A)^{-1}.
\end{equation}
\bigskip

To verify \eqref{eq:may10b24}, observe that
\[
\begin{split}
\frac{1}{2\pi}\int_0^{2\pi}\Delta_{\gE}(t)(R_\omega \wh{F})(e^{it})dt&=\frac{1}{2\pi}\int_0^{2\pi}\Delta_{\gE}(t)\frac{e^{it}F(e^{it})
-\omega F(\omega)}{e^{it}-\omega}dt\\
&=\frac{1}{2\pi}\int_0^{2\pi}\Delta_{\gE}(t)F(e^{it})\frac{\{e^{it}(I_m-\omega A)-\omega(I_m-e^{it}A)\}(I_m-\omega A)^{-1}}
{e^{it}-\omega}dt\\
&=\frac{1}{2\pi}\int_0^{2\pi}\Delta_{\gE}(t)F(e^{it})(I_m-\omega A)^{-1}\\
&={\gee}_0^*G(I_m-\omega A)^{-1}
\end{split}
\]
The verification of \eqref{eq:may10a24} is similar, but simpler.
\bigskip

\noindent
{\bf 3.}\ {\it Verification of \eqref{eq:may9d24}}
\bigskip

Clearly,
\[
F\LL{\gu}_0=2FG{\gu}_0-F\LL^*{\gu}_0=2{\gee}_0^*G{\gu}_0-F\LL^*{\gu}_0.
\]
Thus, as
\[
\begin{split}
F(\omega)\LL^*{\gu}_0&={\gee}_0^*\LL^*(I_m-\omega A)^{-1}{\gu}_0=\{{\gee}_0^*2G-h_0{\gee}_0^*\}(I_m-\omega A)^{-1}{\gu}_0\\
&=2{\gee}_0^*G(I_m-\omega A)^{-1}{\gu}_0-h_0E_+(\omega),
\end{split}
\]
it follows that
\[
\begin{split}
F(\omega)\LL{\gu}_0&=2{\gee}_0^*G{\gu}_0-2{\gee}_0^*G(I_m-\omega A)^{-1}{\gu}_0+h_0E_+(\omega)\\
&=h_0E_+(\omega)+2{\gee}_0^*G\{I_m-(I_m-\omega A)^{-1}\}{\gu}_0\\
&=h_0E_+(\omega)-2\omega{\gee}_0^*GA(I_m-\omega A)^{-1}{\gu}_0,
\end{split}
\]
which justifies the first equality in \eqref{eq:may9d24}. The other equalities are obtained by noting that
\[
2\omega{\gee}_0^*GA(I_m-\omega A)^{-1}{\gu}_0=\frac{2\omega}{2\pi}\int_0^{2\pi}\Delta_{\gE}(t)(R_\omega E_+)(e^{it})dt
\]
and evaluating the last integral as
\[
\frac{2\omega}{2\pi}\int_0^{2\pi}\Delta_{\gE}(t)\frac{E_+(e^{it})-E_+(\omega)}{e^{it}-\omega}dt=\textcircled{1}-\textcircled{2}
\]
with
\[
\textcircled{1}=
\frac{2}{2\pi}\int_0^{2\pi}(E_+(e^{it})^*)^{-1}\frac{\omega e^{-it}}{1-\omega e^{-it}}dt=0
\]
and
\[
\textcircled{2}=\frac{1}{2\pi}\int_0^{2\pi}
\left\{\frac{2\omega}{e^{it}-\omega}+1\right\}\Delta_{\gE}(t)dtE_+(\omega)+h_0E_+(\omega)=
h_0E_+(\omega)-\Phi_{\gE}(\omega)E_+(\omega).
\]

\noindent
{\bf 4.}\ {\it Verification of \eqref{eq:may7a24}}
\bigskip

Clearly,
\[
\begin{split}
\wh{F}(\omega)\LL^*{\gu}_n=&\omega{\gee}_0^*\LL^*(I_m-\omega A)^{-1}{\gu}_n=\omega\{{\gee}_0^*2G-h_0{\gee}_0^*\}(I_m-\omega A)^{-1}{\gu}_n\\
&=2\omega {\gee}_0^*G(I_m-\omega A)^{-1}{\gu}_n-h_0E_-(\omega),
\end{split}
\]
which justifies the first equality in \eqref{eq:may7a24}. The remaining equalities are obtained by 
noting that in view of \eqref{eq:may10b24}, 
\[
2\omega {\gee}_0^*G(I_m-\omega A)^{-1}{\gu}_n=\frac{2\omega}{2\pi}\int_0^{2\pi}\Delta_{\gE}(t)(R_\omega E_-)(e^{it})dt
\]
and then evaluating the integral by Cauchy's formula as
\[
\frac{2\omega}{2\pi}\int_0^{2\pi}\Delta_{\gE}(t)\frac{e^{it}F(e^{it})-\omega F(\omega)}{e^{it}-\omega}dt\,{\gu}_n=	
\textcircled{1}-\textcircled{2},
\]
where
\[
\textcircled{1}=\frac{2\omega}{2\pi}\int_0^{2\pi}\Delta_{\gE}(t)\frac{E_-(e^{it})}{e^{it}-\omega}dt=	\frac{2\omega}{2\pi}\int_0^{2\pi}\frac{({\gu}_n^*F^\#(e^{it}))^{-1}}{e^{it}-\omega}e^{it}dt=2(E_-^\#(\omega))^{-1}
\]
and
\[
\begin{split}
\textcircled{2}&=\frac{1}{2\pi}\int_0^{2\pi}\Delta_{\gE}(t)\left(	\frac{2\omega}{e^{it}-\omega}+1\right)dt\, F(\omega)\,{\gu}_n-\frac{1}{2\pi}\int_0^{2\pi}\Delta_{\gE}(t)dt\, E_-(\omega)\\
&=\Phi_{\gE}(\omega)\,F(\omega)\,{\gu}_n-h_0\,E_-(\omega).
\end{split}
\]
\end{proof}

\begin{rem}
If $(T1)$--$(T3)$ are in force, then the matrix polynomials $E_\pm^\circ$ can also be expressed in terms of the orthogonal projection $\Pi_+$ onto the Hardy space $H_2^p(\DD)$ by the formulas 
\begin{equation}
 \label{eq:nov30bb23}
 E_+^\circ=\Pi_+ \Phi_{\gE}E_+  \quad \textrm{and} \quad E_-^\circ=- \Pi_+\Phi_{\gE}^*E_-.
 \end{equation}
 (In keeping with the convention introduced in \eqref{eq:nov30e23}, interpret $\Pi_+\Delta_{\gE}F$ to mean the projection column by column.)
 \end{rem}

The analysis in \cite{mm03} and \cite{mm05} rests on the observation that the Arov-Grossman model yields 
a one to one correspondence between 
  \begin{equation}
\label{eq:mar19a24}
\left\{\frac{1}{2\pi }\int_{0}^{2\pi} \frac{1}{1-e^{-it}\omega}d\sigma(t):\ \sigma\in\Sigma_T\right\}
\end{equation}
and the set
\begin{equation}
\label{eq:mar19b24}
\{(I-{\cX}-{\cY})^{-1}F{\gee}_0,F{\gee}_0\rangle_{{\cH}_G}:\,
S\in{\cS}^{\ptp}\}
\end{equation}
based on the pair of strictly contractive operators 
$$
{\cX}=\omega R_0\Pi_0\quad\textrm{and}\quad {\cY}=\omega\beta_S(\omega)\Pi_0^\perp,  
$$ 
where 
\begin{equation}
\label{eq:dec8e23}
\beta_S(\omega): \ F{\gu}_0\mapsto -F{\gu}_nS(\omega)\quad\textrm{in which}\ S\in{\cS}^{\ptp}(\DD)\ \textrm{and}\ \omega\in\DD.
\end{equation}
The inverse in \eqref{eq:mar19b24} exists and
\begin{equation}
\label{eq:may1a24}
(I-{\cX}-{\cY})^{-1}Fu=(I-{\cX})^{-1}\{I-{\cY}(I-{\cX})^{-1}\}^{-1}Fu
\end{equation}
because 
${\cX}+{\cY}$ is also a strictly contractive operator:  since 
 the range of ${\cX}$ is orthogonal to the range of ${\cY}$, 
\[
\Vert ({\cX}+{\cY})Fu\Vert_{{\cH}_G}^2=\Vert {\cX}Fu\Vert_{{\cH}_G}^2+\Vert {\cY}Fu\Vert_{{\cH}_G}^2\le \vert\omega\vert^2\left\{ \Vert \Pi_0Fu\Vert_{{\cH}_G}^2+ \Vert \Pi_0^\perp Fu\Vert_{{\cH}_G}^2\right\}=\vert\omega\vert^2 \Vert Fu\Vert_{{\cH}_G}^2.
 \]

Recall that 
 in view of Theorem \ref{thm:feb25a24}, 
the orthogonal projections $\Pi_0$ from ${\cH}_G$ onto $({\cH}_G)_0$ and $\Pi_0^\perp$ from ${\cH}_G$ onto ${\cH}_G\ominus({\cH}_G)_0$ are given by the formulas 
\begin{equation}
\label{eq:dec8b23}
\Pi_0Fu=FN_0Gu \quad and \quad \Pi_0^\perp Fu=F{\gu}_0{\gu}_0^*Gu,
\end{equation}
where 
$N_0=\Gamma-{\gu}_0{\gu}_0^*$ and  $N_\bullet=\Gamma-{\gu}_n{\gu}_n^*$. \begin{lem}
\label{lem:dec8c23}
If $(T1)$--$(T3)$ are in force, then
\begin{equation}
\label{eq:may8b24}
\Vert AN_0G\Vert=\Vert N_0G\Vert\ge 1, \quad \Vert G^{1/2}N_0G^{1/2}\Vert= 1 \quad\textrm{and}\quad \sum_{k=0}^\infty \Vert {\gu}_0^*G(AN_0G)^kA{\gu}_0\Vert<\infty.
\end{equation}
Moreover,
\begin{equation}
\label{eq:dec8d23}
E_+(\omega)^{-1}=\{I_p-\omega {\gu}_0^*G A(I_m-\omega N_0GA)^{-1}{\gu}_0\} (\gamma_{00}^{(n)})^{-1/2}
\end{equation}
and
\begin{equation}
\label{eq:may5c24}
{\gu}_0E_+(\omega)^{-1}=(I_m-\omega A)(I_m-\omega N_0GA)^{-1}{\gu}_0(\gamma_{00}^{(n)})^{-1/2}
\end{equation}
for every point $\omega$ in the closed unit disk $\ol{\DD}$.
\end{lem}

\begin{proof}
Since $I_m={\gee}_0{\gee}_0^*+A^*A$ and ${\gee}_0^*N_0=0$, it is readily seen that 
\[
\Vert N_0G\Vert=\Vert A^*AN_0G\Vert\le \Vert AN_0G\Vert\le\Vert N_0G\Vert.
\]
Therefore, $\Vert AN_0G\Vert=\Vert N_0G\Vert$. Moreover, $\Vert N_0G\Vert\ge 1$, since $\Vert N_0G\Vert=\Vert (N_0G)^2\Vert\le \Vert N_0G\Vert^2$ and $N_0G\ne O$. 
 The equality 
$ \Vert G^{1/2}N_0G^{1/2}\Vert= 1 $ is justified in Theorem \ref{thm:may14a24}.

By definition,
\[
\begin{split} 
E_+(\omega)&={\gee}_0^*(I_m-\omega A)^{-1}{\gu}_0={\gee}_0\{(I_m-\omega A)^{-1}-I_m\}{\gu}_0+{\gee}_0^*{\gu}_0\\
&={\gee}_0^*(I_m-\omega A)^{-1}\omega A{\gu}_0+ (\gamma_{00}^{(n)})^{1/2}\\
&=\{I_p+\omega{\gee}_0^*A(I_m-\omega A)^{-1}{\gu}_0(\gamma_{00}^{(n)})^{-1/2}
\}(\gamma_{00}^{(n)})^{1/2}.
\end{split}
\]
The term in curly brackets is of the form $I_p+\omega C(I_m-\omega A)^{-1}B$ with $C={\gee}_0^*A$ and $B={\gu}_0(\gamma_{00}^{(n)})^{-1/2}$. Thus, as is well known and is readily checked,
\[
\{I_p+\omega C(I_m-\omega A)^{-1}B\}\{I_p-\omega C(I_m-\omega A^\times)^{-1}B\}=I_p\quad\textrm{when}\quad 
A^\times=A-BC.
\]
But this yields  \eqref{eq:dec8d23}, since $A-BC=N_0GA$ and ${\gu}_0^*G=(\gamma_{00}^{(n)})^{-1/2}{\gee}_0^*$.
Moreover, as \eqref{eq:dec8d23} can also be expressed as
$$
E_+(\omega)^{-1}=\{I_p-\omega {\gu}_0^*G (I_m-\omega AN_0G)^{-1}{\gu}_0\}A (\gamma_{00}^{(n)})^{-1/2}
$$
and $\det\,E_+(\omega)\ne 0$ in the closed unit disc, a theorem of Wiener guarantees that $E_+^{-1}$ belongs to the Wiener algebra of $\ptp$ mvf's on the the unit circle with summable Fourier coefficients. Therefore, the second bound in 
\eqref{eq:may8b24} holds.

Finally, in view of \eqref{eq:dec8d23}, 
\[
\begin{split}
{\gu}_0E_+(\omega)^{-1}&={\gu}_0\{I_p-\omega {\gu}_0^*G A(I_m-\omega N_0GA)^{-1}{\gu}_0\} (\gamma_{00}^{(n)})^{-1/2}\\
&=\{I_m-\omega{\gu}_0{\gu}_0^*GA(I_m-\omega N_0GA)^{-1}\}{\gu}_0(\gamma_{00}^{(n)})^{-1/2}\\
&=\{I_m-\omega N_0GA-\omega{\gu}_0{\gu}_0^*GA\}(I_m-\omega N_0GA)^{-1}{\gu}_0(\gamma_{00}^{(n)})^{-1/2},
\end{split}
\]
which coincides with \eqref{eq:may5c24}, since
\[
I_m-\omega N_0GA-\omega{\gu}_0{\gu}_0^*GA= I_m-\omega(\Gamma-{\gu}_0{\gu}_0^*)GA-
\omega{\gu}_0{\gu}_0^*GA=I_m-\omega A.
\]
\end{proof}

\begin{lem} 
\label{lem:dec7c23}
If $(T1)$--$(T3)$ are in force, then the following identities are valid for every choice of $\omega\in\DD$ and $u\in\CC^m$:

\begin{equation}
\label{eq:dec7f23}
(I-\omega R_0\Pi_0)^{-1}Fu=F(I_m-\omega AN_0G)^{-1}u
\end{equation}

\begin{equation}
\label{eq:dec8c23}
\Pi_0^\perp(I-\omega R_0\Pi_0)^{-1}Fu=E_+ E_+(\omega)^{-1}F(\omega)u 
\end{equation}

\end{lem}

\begin{proof}
To evaluate \eqref{eq:dec7f23}, 
let
$(I-\omega R_0\Pi_0)^{-1}Fu=Fv$. Then 
\[
\begin{split}
Fu&=(I-\omega R_0\Pi_0)Fv=Fv-\omega R_0FN_0Gv\\
&=F\{I_m-\omega AN_0G\}v.
\end{split}
\]
Therefore, $u=(I_m-\omega A N_0G)v$, and hence $v=(I_m-\omega AN_0G)^{-1}v$, 
which leads easily to  \eqref{eq:dec7f23}.

Next, to evaluate \eqref{eq:dec8c23}, observe that
in view of \eqref{eq:dec7f23} and \eqref{eq:dec8b23},
\[
\begin{split}
\Pi_0^\perp(I-\omega R_0\Pi_0)^{-1}Fu&=F{\gu}_0{\gu}_0^*G(I_m-\omega AN_0G)^{-1}u\\
&=F{\gu}_0{\gu}_0^*G(I_m-\omega A+\omega A{\gu}_0{\gu}_0^*G)^{-1}u\\
&=E_+{\gu}_0^*G(I_m+(I_m-\omega A)^{-1}\omega A{\gu}_0{\gu}_0^*G)^{-1}(I_m-\omega A)^{-1}u\\
&=E_+(I_p+{\gu}_0^*G(I_m-\omega A)^{-1}\omega A{\gu}_0)^{-1}{\gu}_0^*G(I_m-\omega A)^{-1}u\\
&=E_+({\gu}_0^*G[I_m+(I_m-\omega A)^{-1}\omega A]{\gu}_0)^{-1}{\gu}_0^*G(I_m-\omega A)^{-1}u\\
&=E_+({\gu}_0^*G(I_m-\omega A)^{-1}{\gu}_0)^{-1}{\gu}_0^*G(I_m-\omega A)^{-1}u\\
&=E_+E_+(\omega)^{-1}F(\omega)u,
\end{split}
\]
since ${\gu}_0^*G=(\gamma_{00}^{(n)})^{-1/2}{\gee}_0^*$.
\end{proof}

\begin{lem}
\label{lem:may8b24}
If $(T1)$--$(T3)$ are in force,  then 
\begin{equation}
\label{eq:dec8f23}
\begin{split}
\{(I-{\cY}(I-{\cX})^{-1}\}^{-1}Fu&=
\{I-\omega\beta_S(\omega)\Pi_0^\perp(I-\omega R_0\Pi_0)^{-1}\}^{-1}Fu\\ &=Fu-F{\gu}_n\{\omega S(\omega)(E_+(\omega)+ E_-(\omega)  S(\omega))^{-1}F(\omega)\}u
\end{split}
\end{equation}
for every choice of $u\in\CC^m$ and $\omega\in\DD$. 
\end{lem}

\begin{proof}
Let
\[
\{I-{\cY}(I-{\cX})^{-1}\}^{-1}Fu=Fv.
\]
Then, by \eqref{eq:dec8c23}, 
\[
\begin{split}
Fu&=
\{I-{\cY}(I-{\cX})^{-1}\}Fv
=Fv-\omega \beta_S(\omega)E_+E_+(\omega)^{-1}F(\omega)v\\ &=Fv+\omega F{\gu}_nS(\omega)E_+(\omega)^{-1}F(\omega)v
=F\{I_m+\omega{\gu}_nS(\omega)E_+(\omega)^{-1}F(\omega)\}v\\
&=F\{I_m-BC\}v 
\end{split}
\]
with
\[
B=-\omega {\gu}_nS(\omega)
\quad\textrm{ and}\quad  C=E_+(\omega)^{-1}F(\omega).
\]
Thus, as 
$$
\Vert CB\Vert=\Vert E_+(\omega)^{-1}\omega F(\omega){\gu}_nS(\omega)\Vert\le
\Vert E_+(\omega)^{-1} E_-(\omega)\Vert<1\quad\textrm{for}\ \omega\in\DD,
$$
the matrices $I_p-CB$ and $I_m-BC$ are both invertible. Moreover, as 
\[
(I_m-BC)^{-1}=I_m+B(I_p-CB)^{-1}C, 
\]
the condition $Fu=F(I_m-BC)v$ implies that $u=(I_m-BC)v$ and hence that 
\[
\begin{split}
v&=(I_m-BC)^{-1}u=(I_m+B(I_p-CB)^{-1}C)u\\
&=\{I_m-\omega {\gu}_nS(\omega)(E_+(\omega)+  E_-(\omega)  S(\omega))^{-1}F(\omega)\}u,
\end{split}
\]
which serves to justify \eqref{eq:dec8f23}.
\end{proof}

\begin{lem}
\label{lem:may9a24}
If $(T1)$--$(T3)$ are in force,  then 
\begin{equation}
\label{eq:mar19aa24}
{\gee}_0^* G(I_m-\omega AN_0G)^{-1}{\gu}_n=\{\omega E_-^\#(\omega)\}^{-1} 
\end{equation}
\end{lem}

\begin{proof}
The formula in  (2) of Corollary \ref{cor:mar20a24} is equivalent to  the identity 
\[
\omega AN_0F^{\#}(\omega)=N_\bullet F^{\#}(\omega).
\]
Consequently, 
\[
\begin{split}
(I_m-\omega AN_0G)^{-1}{\gu}_n{\gu}_n^*F^\#(\omega)&=(I_m-\omega AN_0G)^{-1}(\Gamma-N_\bullet)F^\#(\omega)\\
&=\Gamma (\Gamma-\omega AN_0)^{-1}(\Gamma-\omega AN_0)F^\#(\omega)=\Gamma F^\#(\omega).
\end{split}
\]
Therefore,
\[
{\gee}_0^*G(I_m-\omega AN_0G)^{-1}{\gu}_n{\gu}_n^*F^\#(\omega)={\gee}_o^*G\Gamma F^\#(\omega)=I_p.
\]
Thus, 
\[
{\gee}_0^*G(I_m-\omega AN_0G)^{-1}{\gu}_n=\{{\gu}_n^*F^\#(\omega)\}^{-1}=\{{\gu}_n^*\omega\wh{F}^\#(\omega)\}^{-1}
=\{\omega E_-^\#(\omega)\}^{-1}
\]
since  $E_-^\#={\gu}_n^*\wh{F}^\#$.
\end{proof}

\begin{lem}
\label{lem:may9b24}
If $(T1)$--$(T3)$ are in force,  then 
\begin{equation}
\label{eq:may5a24}
\Phi_{\gE}(\omega)=2{\gee}_0^*G(I_m-\omega A N_0G)^{-1}{\gee}_0-h_0={\gee}_0^*G(I_m+\omega AN_0G)(I_m-\omega AN_0G)^{-1}{\gee}_0.
\end{equation}
\end{lem}

\begin{proof}
The asserted verification rests on the identity
$\Phi_{\gE}(\omega)=E_+^\circ(\omega)E_+(\omega)^{-1}$,  the first equality in \eqref{eq:may9d24}.
and formula \eqref{eq:dec8d23}:
\[
\begin{split}
E_+^\circ(\omega)E_+(\omega)^{-1} &=h_0-2\omega {\gee}_0^*GA(I_m-\omega A)^{-1}{\gu}_0 E_+(\omega)^{-1}\\
&=h_0-2\omega {\gee}_0^*GA(I_m-\omega A)^{-1}(I_m-\omega A)(I_m-\omega N_0GA)^{-1}{\gu}_0(\gamma_{00}^{(n)})^{-1/2}\\
&=h_0-2\omega {\gee}_0^*GA(I_m-\omega N_0GA)^{-1}{\gu}_0(\gamma_{00}^{(n)})^{-1/2}\\
&=h_0-2 {\gee}_0^*G(I_m-\omega A N_0G)^{-1}\omega A{\gu}_0(\gamma_{00}^{(n)})^{-1/2}.
\end{split}
\]
But this in turn implies that 
\begin{equation}
\label{eq:may6a24}
\Phi_{\gE}(\omega)=h_0+2{\gee}_0^*G(I_m-\omega A N_0G)^{-1}\omega AN_0G{\gee}_0,
\end{equation}
since
\[
A{\gu}_0(\gamma_{00}^{(n)})^{-1/2}=A{\gu}_0(\gamma_{00}^{(n)})^{-1/2}{\gee}_0^*\Gamma G{\gee}_0=A{\gu}_0{\gu}_0^*G{\gee}_0=-AN_0G{\gee}_0.
\]
Thus, as \eqref{eq:may6a24} is equivalent to \eqref{eq:may5a24}, the proof is complete.
\end{proof}

\begin{lem}
\label{lem:may12b24}
If $T(1)$--$T(3)$ are in force, then
\begin{equation}
\label{eq:mar20b24}
\begin{split}
{\gee}_0^*G(I_m-\omega AN_0G)^{-1}{\gee}_0&=\frac{1}{2\pi}\int_0^{2\pi}\frac{1}{1-\omega e^{-it}}\Delta_{\gE}(t)dt\\
&=\frac{\Phi_{\gE}(\omega)+\Phi_{\gE}(0)}{2} \quad for\ \omega\in\DD
\end{split}
\end{equation}
and hence
\begin{equation}
\label{eq:may5e24}
\frac{1}{2\pi}\int_0^{2\pi}e^{-itk}\Delta_{\gE}(t)dt={\gee}_0^*G(AN_0G)^k{\gee}_0\quad for\ k=0,1,\ldots .
\end{equation}
\end{lem}

\begin{proof}
Formula \eqref{eq:mar20b24} drops out by comparing \eqref{eq:may5a24} with the integral representation 
\[
\begin{split}
\Phi_{\gE}(\omega)&=\frac{1}{2\pi}\int_0^{2\pi}[2-(1-\omega e^{-it})](1-\omega e^{-it})^{-1}\Delta_{\gE}(t)dt\\
&=\frac{2}{2\pi}\int_0^{2\pi}(1-\omega e^{-it})^{-1}\Delta_{\gE}(t)dt-h_0.
\end{split}
\]
Formula \eqref{eq:may5e24} is immediate from \eqref{eq:mar20b24}. 
\end{proof}

\begin{lem}
\label{lem:may8c24}
If $(T1)$--$(T3)$ are in force,  then 
\begin{equation}
\label{eq:may1b24}
\begin{split}
\langle (I-{\cX}-{\cY})^{-1}F{\gee}_0 x, F{\gee}_0 y\rangle_{{\cH}_G}&=
\{\Phi_{\gE}(\omega)+h_0\}/2\\ &-y^*E_-^{\#}(\omega)^{-1}S(\omega)\{E_-(\omega)S(\omega)+E_+(\omega)\}^{-1}x
\end{split}
\end{equation}
for every choice of $x,y\in\CC^p$ and $\omega\in\DD$.
\end{lem}

\begin{proof}
In view of \eqref{eq:may1a24} and \eqref{eq:dec8f23}, 
\[
(I-{\cX}-{\cY})^{-1}Fu=(I-{\cX})^{-1}\{Fu-F{\gu}_n\omega Y(\omega)F(\omega)\}u
\]
with $$Y(\omega)=S(\omega)(E_-(\omega)S(\omega)+E_+(\omega))^{-1}$$ for short. Thus,
 the left hand side of \eqref{eq:may1b24} is equal to $\textcircled{1}-\textcircled{2}$, where
\[
\textcircled{1}=\langle (I-{\cX})^{-1}F{\gee}_0x,F{\gee}_0y\rangle_{{\cH}_G}=y^*{\gee}_0^*G(I_m-\omega AN_0G)^{-1}{\gee}_0 x=y^*\left\{\frac{\Phi_{\gE}(\omega)+h_0}{2}\right\}x
\]
and
\[
\textcircled{2}=\langle (I-{\cX})^{-1}F{\gu}_n\omega Y(\omega)F(\omega){\gee}_0x,F{\gee}_0y\rangle_{{\cH}_G}=y^*{\gee}_0^*G(I_m-\omega AN_0G)^{-1}{\gu}_n\omega Y(\omega)x
\]
(since $F(\omega){\gee}_0x=x$).
\end{proof}

\begin{thm}
\label{thm:may7a24}
If $(T1)$--$(T3)$ are in force,
 then the formula
\begin{equation}
\left\{\frac{1}{2\pi }\int_{0}^{2\pi} \frac{e^{it}+\omega}{e^{it}-\omega}d\sigma(t):\ \sigma\in\Sigma_{\rm T} \right\}
=\{(E_-^\circ S+E_+^\circ)(E_-S+E_+)^{-1}:\ S\in{\cS}^{\ptp}(\DD)\}
\end{equation}
establishes a $1:1$ correspondence between the solutions $\sigma$ of the truncated matrix trigonometric moment problem and the mvf's $S$ in the Schur class ${\cS}^{\ptp}(\DD)$.
\end{thm}

\begin{proof} The proof rests on the fact that the set considered in \eqref{eq:mar19a24} coincides with the set considered in \eqref{eq:mar19b24} .
Let ${\cZ}={\cX}+{\cY}$. Then, in terms of the notation in the proof of Lemma \ref{lem:may8c24},
\[
\begin{split}
\langle (I+{\cZ})(I-{\cZ})^{-1}F{\gee}_0,F{\gee}_0\rangle_{{\cH}_G}&
=2\langle (I-{\cZ})^{-1}F{\gee}_0,F{\gee}_0\rangle_{{\cH}_G}-
\langle F{\gee}_0,F{\gee}_0\rangle_{{\cH}_G}\\
&=2\{\textcircled{1}-\textcircled{2}\}-\langle F{\gee}_0,F{\gee}_0\rangle_{{\cH}_G}\\
&=\Phi_{\gE}-2(E_-^\#)^{-1}S(E_-S+E_+)^{-1}\\
&=\{\Phi_{\gE}(E_-S+E_+)-2(E_-^\#)^{-1}S\}(E_-S+E_+)^{-1}\\
&=(E_-^\circ S+E_+)(E_-S+E_+)^{-1}.
\end{split}
\]
\end{proof}

\section{Solutions to the truncated matrix Hamburger moment problem}
\label{sec:solutionshamburger}

Recall that when $G$ is block Hankel, then ${\cH}_G$ is a de Branges space ${\cB}({\gE})$ with ${\gE}=\begin{bmatrix}E_-&E_+\end{bmatrix}$,
\begin{equation}
\label{eq:dec1b23}
E_+ (\lambda)= \frac{\rho_\alpha(\lambda)}{\sqrt{\rho_\alpha(\alpha)}}F(\lambda){\gz}_\alpha, \quad 
		E_ - (\lambda)= \frac{\rho_{\ol{\alpha}}(\lambda)}{\sqrt{\rho_\alpha(\alpha)}}F(\lambda){\gz}_{\ol{\alpha}}, \quad
\rho_\omega(\lambda)=-2\pi i(\lambda-\ol{\omega})
\end{equation}
 and inner product \eqref{eq:dec12g23}.

 The main objective of this section is to develop a
 description of the set 
\[
\Sigma_{\rm H}\quad \textrm{of solutions $\sigma$ to the truncated matrix Hamburger moment problem}
\]
analogous to Theorem \ref{thm:may7a24}
 in terms  of the matrix polynomials $E_\pm$  defined in \eqref{eq:dec1b23} and the mvf's (which, as we shall see shortly,  are also matrix polynomials) 
\begin{equation}
\label{eq:may21a24}
E_\pm^\circ(\omega)=-\frac{1}{\pi i}\int_{-\infty}^\infty \Delta_{\gE}(\mu) (R_\omega E_\pm)(\mu)d\mu
\end{equation}
The proof again rests on the Arov-Grossman parameterization of minimal unitary extensions of isometric operators applied to the isometric operator $T_\alpha=I+(\alpha-\ol{\alpha})R_\alpha$ considered in Theorem \ref{thm:sep23b20} and a theorem of Chumakin \cite{ch67} on generalized resolvents of isometric operators. 

The analysis, which is similar in spirit to the operator theory approach presented in \cite{za12}, rests on Theorem \ref{thm:jun4a24}, which is formulated below in terms of the unitary matrix $\chi_\infty=\lim_{\nu\uparrow\infty} (E_+^{-1}E_-)(i\nu)$,  the set  

\begin{equation}
\label{eq:jun4a24}
{\cS}_{\gE}^{\ptp}(\CC_+)=\{S\in{\cS}^{\ptp}(\CC_+):\ \lim_{\nu\uparrow\infty} \nu^{-1}\{I_p+\chi_\infty S(i\nu)\}^{-1}=0\}
\end{equation}
and a pair of strictly contractive operators ${\cX}=b_\alpha(\omega)T_\alpha\Pi_\alpha$ and ${\cY}=b_\alpha(\omega)\beta_S(\omega)\Pi_\alpha^\perp$ that are defined in terms of the elementary Blaschke factor 
$b_\alpha(\lambda)=(\lambda-\alpha)/(\lambda-\ol{\alpha})$ with $\alpha\in\CC_+$ and the operator 
$$
\beta_S(\omega): \ F{\gz}_\alpha\mapsto -F{\gz}_{\ol{\alpha}}S(\omega)\quad\textrm{in which}\ S\in{\cS}^{\ptp}(\CC_+)\ \textrm{and}\ \omega\in\CC_+.
$$
However, once again we use the underlying de Branges space structure to compute projections in terms of appropriately chosen reproducing kernels rather than orthogonal bases. This eases the bookkeeping and yields explicit formulas for the coefficients in the linear fractional transformation that serves to parameterize the solutions of the truncated matrix Hamburger moment problem. Our starting point is:

\begin{thm}
\label{thm:jun4a24}
The formula 
\begin{equation}
\label{eq:jun4b24}
\frac{1}{\pi i}\int_{-\infty}^\infty\frac{d\sigma(\mu)}{\mu-\omega}=\frac{2\rho_\alpha(\alpha)}{\rho_\alpha(\omega)\rho_{\ol{\alpha}}(\omega)}\left\{\frac{\rho_\alpha(\omega)}{\rho_\alpha(\alpha)}h_0-\langle(I-{\cX}-{\cY})^{-1}F{\gee}_0,F{\gee}_0\rangle_{{\cH}_G}\right\},\quad \alpha,\omega\in\CC_+,
\end{equation}
establishes a $1:1$ correspondence between the set of solutions $\sigma$ of the truncated matrix Hamburger moment problem and the set of 
$S$ in the restricted Schur class ${\cS}_{\gE}^{\ptp}(\CC_+)$  (${\cY}$ depends upon $S$).
\end{thm}

Formula \eqref{eq:jun4b24} is the counterpart of formula (38) in \cite{za12} in the present setting. 
The operator $I-{\cX}-{\cY}$ that intervenes in \eqref{eq:jun4b24} is invertible because 
${\cX}+{\cY}$ is also a strictly contractive operator in ${\cH}_G$:  
since the range of ${\cX}$ is orthogonal to the range of ${\cY}$,
\[
\Vert ({\cX}+{\cY})Fu\Vert_{{\cH}_G}^2=\Vert {\cX}Fu\Vert_{{\cH}_G}^2+\Vert {\cY}Fu\Vert_{{\cH}_G}^2\le \vert b_\alpha(\omega)\vert^2\left\{ \Vert \Pi_\alpha Fu\Vert_{{\cH}_G}^2+ \Vert \Pi_\alpha^\perp Fu\Vert_{{\cH}_G}^2\right\}=\vert b_\alpha(\omega)\vert^2 \Vert Fu\Vert_{{\cH}_G}^2.
 \]

The formulas
\[
\Pi_\alpha Fu=FN_\alpha Gu,\quad \Pi_\alpha^\perp Fu=F{\gz}_\alpha{\gz}_\alpha^*Gu,\quad \beta_S(\omega)F{\gz}_\alpha=-F{\gz}_{\ol{\alpha}}S(\omega)
\]
and
\begin{equation}
\label{eq:dec1c23}
\Phi_{\gE}(\omega)=\frac{1}{\pi i}\int_{-\infty}^{\infty} \frac{\Delta_{\gE}(\mu)}{\mu-\omega}d\mu\quad\textrm{where}\quad 
\Delta_{\gE}(\mu)=(E_+(\mu)E_+(\mu)^*)^{-1}
\end{equation}
will play a prominent role in the analysis.

For ease of future reference we list the assumptions and definitions that will be in force throughout this section and in future discussion of the truncated matrix Hamburger moment problem.
\begin{enumerate}
\item[\rm(H1)] $G$ is a positive definite block Hankel matrix and $\Gamma=G^{-1}$.
\item[\rm(H2)] $\rho_\omega(\lambda)=-2\pi i(\lambda-\ol{\omega})$ and $f^\#(\lambda)=f(\ol{\lambda})^*$
\item[\rm(H3)] The mvf's $E_\pm$, $E_\pm^\circ$ and $\Phi_{\gE}$   are specified by formulas \eqref{eq:dec1b23}, 
\eqref{eq:may21a24} and \eqref{eq:dec1c23}, respectively, and the inner product in ${\cB}(\gE)$ is given by \eqref{eq:dec12g23} (and \eqref{eq:aug24b23}).
\end{enumerate}

The first order of business is to develop other formulas for the mvf's $E_\pm^\circ$. 
\begin{lem}
\label{lem:nov12a23}
If the positive definite matrix $G$ is block Hankel and if $F(\omega)=\gee_0^\ast (I-\omega A)^{-1}$ and 
	$\wh{F}(\omega)=\omega F(\omega)$, then:
\begin{equation}
\label{eq:nov13e23}
\int_{-\infty}^\infty \Delta_{\gE}(\mu) (R_\omega F)(\mu) d\mu
=  \gee_0^\ast GA (I_m-\omega A)^{-1}
\end{equation}
and
\begin{equation}
\label{eq:nov12a23}
\int_{-\infty}^\infty \Delta_{\gE}(\mu) (R_\omega \wh{F})(\mu) d\mu
= \gee_0^\ast G (I_m-\omega A)^{-1}.
\end{equation}
\end{lem}

\begin{proof}
Under the given assumptions,
\begin{align*}
	 \int \limits_{-\infty}^{\infty} \Delta_{\gE} (\mu )(R_\omega F)(\mu)  d\mu  
		&=  \int \limits_{-\infty}^{\infty} \Delta_{\gE} (\mu ) \gee_0^\ast (I_m-\mu A)^{-1}A(I_m-\omega A)^{-1} d\mu \\
		&= \gee_0^\ast GA (I_m-\omega A)^{-1},
			\end{align*}
whereas,
\begin{align*}
	  \int_{-\infty}^{\infty} \Delta_{\gE} (\mu )(R_\omega \wh{F})(\mu) d\mu 
&= \int_{-\infty}^{\infty} \Delta_{\gE} (\mu )\{F(\mu)+\omega(R_\omega F)(\mu) \} d\mu \\
&={\gee}_0^*G+\omega  \gee_0^\ast GA (I_m-\omega A)^{-1}=  \gee_0^\ast G (I_m-\omega A)^{-1},		
\end{align*}
as claimed. \end{proof}

\begin{lem}
\label{lem:nov19a23}

The mvf's $E_\pm^\circ$  are matrix polynomials (that are often referred to as matrix polynomials of the second kind) that can also be expressed as 
\begin{equation}
\label{eq:nov19b23}
E_+^\circ(\omega)=\frac{2}{\sqrt{\rho_\alpha(\alpha)}}{\gee}_0^*G(I_m-\omega A)^{-1}(I_m-\ol{\alpha} A){\gz}_\alpha
=\Phi_{\gE}(\omega)E_+(\omega)
\end{equation}
and
\begin{equation}
\label{eq:may20d24}
E_-^\circ(\omega)=\frac{2}{\sqrt{\rho_\alpha(\alpha)}}{\gee}_0^*G(I_m-\omega A)^{-1}(I_m-\alpha A){\gz}_{\ol{\alpha}}.
=\Phi_{\gE}(\omega)E_-(\omega)-2(E_-^\#(\omega))^{-1}
\end{equation}
when $\omega\in\CC_+$. 
\end{lem}

\begin{proof} 
The first equalities in \eqref{eq:nov19b23} and \eqref{eq:may20d24} follow easily from Lemma \ref{lem:nov12a23}, since
\[
E_+=-\frac{2\pi i}{\sqrt{\rho_\alpha(\alpha)}}(\wh{F}-\ol{\alpha}F){\gz}_\alpha\quad\textrm{and}\quad 
E_-=-\frac{2\pi i}{\sqrt{\rho_\alpha(\alpha)}}(\wh{F}-\alpha F){\gz}_{\ol{\alpha}}.
\]
These formulas serve also to show that $E_\pm^\circ$ are matrix polynomials. The second set of equalities is obtained by exploiting the fact that the vvf's $E_+^{-1}\xi, (E_-^\#)^{-1}\xi$ belong to $H_2^p(\CC_+)$ for every $\xi\in\CC^p$ in order to evaluate the integrals in \eqref{eq:may21a24} by Cauchy's formula: 
 In view of formulas \eqref{eq:dec1c23},  \eqref{eq:nov15e23} and \eqref{eq:nov12a23}, 
\[
\begin{split}
E_+^\circ(\omega)&=-\frac{1}{\pi i}\int_{-\infty}^\infty \Delta_{\gE}(\mu) (R_\omega E_+)(\mu)d\mu=
-\frac{1}{\pi i}\int_{-\infty}^\infty \frac{\Delta_{\gE}(\mu)}{\mu-\omega}E_+(\mu)d\mu+\Phi_{\gE}(\omega)E_+(\omega)\\&=\Phi_{\gE}(\omega)E_+(\omega)\quad\textrm{for}\ \omega\in\CC_+,
\end{split}
\]
since $\Delta_{\gE}(\mu)=\{E_+(\mu)E_+(\mu)^*\}^{-1}=\{E_-(\mu)E_-(\mu)^*\}^{-1}$ for  $\mu\in\RR$ and 
$$ 
\frac{1}{\pi i}\int_{-\infty}^\infty \frac{\Delta_{\gE}(\mu)}{\mu-\omega}E_+(\mu)d\mu=\frac{1}{\pi i}\int_{-\infty}^\infty (E_+(\mu)^*)^{-1}\frac{1}{\mu-\omega}d\mu=0 \quad\textrm{when}\ \omega\in\CC_+.
$$
Analogously,
$$
E_-^\circ(\omega)=-\frac{1}{\pi i}\int_{-\infty}^\infty \frac{\Delta_{\gE}(\mu)}{\mu-\omega}E_-(\mu)d\mu+\Phi_{\gE}(\omega)E_-(\omega)=-2(E_-^\#(\omega))^{-1}+\Phi_{\gE}(\omega)E_-(\omega)\quad\textrm{for}\ \omega\in\CC_+,
$$
since
$$ 
\frac{1}{\pi i}\int_{-\infty}^\infty \frac{\Delta_{\gE}(\mu)}{\mu-\omega}E_-(\mu)d\mu=\frac{1}{\pi i}\int_{-\infty}^\infty (E_-(\mu)^*)^{-1}\frac{1}{\mu-\omega}d\mu=2E_-^\#(\omega)^{-1} \quad\textrm{when}\ \omega\in\CC_+.
$$
\end{proof}

\begin{rem}
\label{rem:may21a24}
The formulas in \eqref{eq:may21a24} can also be expressed in terms of the orthogonal projection $\Pi_+$ onto $H_2^p(\CC_+)$, since
\begin{equation}
\label{eq:nov15d23}
\begin{split}
\Pi_+\Phi_{\gE}^* Fu&=2\Pi_+ \Delta_{\gE}F-\Phi_{\gE}F=\frac{2}{2\pi i}\int_{-\infty}^\infty\frac{\Delta_{\gE}(\mu)}{\mu-\omega}F(\mu)d\mu-
\frac{2}{2\pi i}\int_{-\infty}^\infty\frac{\Delta_{\gE}(\mu)}{\mu-\omega}d\mu F(\omega)\\ &=
\frac{1}{\pi i}\int_{-\infty}^\infty \Delta_{\gE}(\mu) (R_\omega F)(\mu) d\mu.
\end{split}
\end{equation}
    \end{rem}

\begin{lem}
\label{lem:may20a24}
If $(H1)$--$(H3)$ are in force, then $E_+(\lambda)$ is invertible for every point $\lambda$ in the closed upper half-plane $\ol{\CC_+}$ and   (hence) the matrix ${\gee}_0^*{\gz}_\alpha$ is invertible. Moreover, if $\det\,E_+(\lambda)\ne 0$, then 
\begin{equation}
\label{eq:may20a24}
E_+(\lambda)^{-1}=\frac{\sqrt{\rho_\alpha(\alpha)}}{\rho_\alpha(\lambda)}({\gee}_0^*{\gz}_\alpha)^{-1}\{I_p-\lambda{\gee}_0^*A(I_m-\lambda QA)^{-1}{\gz}_\alpha({\gee}_0^*{\gz}_\alpha)^{-1}\},
\end{equation}
\begin{equation}
\label{eq:may20b24}
{\gz}_\alpha E_+(\lambda)^{-1}=\frac{\sqrt{\rho_\alpha(\alpha)}}{\rho_\alpha(\lambda)}(I_m-\lambda A)(I_m-\lambda QA)^{-1}{\gz}_\alpha({\gee}_0^*{\gz}_\alpha)^{-1}
\end{equation}
and
\begin{equation}
\label{eq:may21f24} 
Q=I_m-{\gz}_\alpha({\gee}_0^*{\gz}_\alpha)^{-1}{\gee}_0^*=Q^2.
\end{equation}
\end{lem}

\begin{proof}
The formula
\begin{equation}
\label{eq:jun14a24}
E_+(\omega)E_+(\omega)^*-E_-(\omega)E_-(\omega)^*=\rho_\omega(\omega)F(\omega)\Gamma F(\omega)^*
\end{equation}
guarantees that $E_+(\omega)$ is invertible for every point $\omega\in \ol{\CC_+}$, since the RK $K_\omega(\omega)=F(\omega)\Gamma F(\omega)^*$ is positive definite for every $\omega\in\CC$. (If  $\omega\in\RR$, then the asserted invertibility  follows from the fact that \eqref{eq:jun14a24} $\implies E_+(\omega)^*x=0\iff E_-(\omega)^*x=0$.) 
Thus, 
${\gee}_0^*{\gz}_\alpha=F(0){\gz}_\alpha=(\sqrt{\rho_\alpha(\alpha)}/\rho_\alpha(0))E_+(0)$ 
is invertible. Consequently,
\[
\begin{split}
E_+(\lambda)
&=\frac{\rho_\alpha(\lambda)}{\sqrt{\rho_\alpha(\alpha)}}\{(F(\lambda)-F(0)){\gz}_\alpha+F(0){\gz}_\alpha\}\\
&=\frac{\rho_\alpha(\lambda)}{\sqrt{\rho_\alpha(\alpha)}}\{\lambda F(\lambda)A{\gz}_\alpha+F(0){\gz}_\alpha\}\\ &=\frac{\rho_\alpha(\lambda)}{\sqrt{\rho_\alpha(\alpha)}}F(0){\gz}_\alpha \{I_p+(F(0){\gz}_\alpha)^{-1}\lambda {\gee}_0^*(I_m-\lambda A)^{-1}A{\gz}_\alpha\},
\end{split}
\]
and the term in curly brackets $\{\cdots\}$ is of the form $I_p+\lambda C(I_m-\lambda A)^{-1}B$ with 
$C=(F(0){\gz}_\alpha)^{-1}{\gee}_0^*$ and $B=A{\gz}_\alpha$. Therefore (much as in the proof of Lemma \ref{lem:dec8c23}),
\[
\{\cdots\}^{-1}=I_p-\lambda C(I_m-\lambda A^\times)^{-1} B
\quad\textrm{with}\ A^\times =A-BC=A(I_m-{\gz}_\alpha ({\gee}_0^*{\gz}_\alpha)^{-1}{\gee}_0^*)=AQ.
\]
Formula \eqref{eq:may20a24} drops out upon combining terms.

Formula \eqref{eq:may20b24} is obtained from \eqref{eq:may20a24} by noting that
\[
\begin{split}
{\gz}_\alpha({\gee}_0^*{\gz}_\alpha)^{-1}\{I_p-&\lambda{\gee}_0^*A(I_m-\lambda QA)^{-1}
{\gz}_\alpha({\gee}_0^*{\gz}_\alpha)^{-1}\}\\
&=\{I_m-\lambda {\gz}_\alpha({\gee}_0^*{\gz}_\alpha)^{-1}{\gee}_0^*A(I_m-\lambda QA)^{-1}\}
{\gz}_\alpha({\gee}_0^*{\gz}_\alpha)^{-1}\\
&=\{I_m-\lambda QA-\lambda {\gz}_\alpha({\gee}_0^*{\gz}_\alpha)^{-1}{\gee}_0^*A\}(I_m-\lambda QA)^{-1}{\gz}_\alpha({\gee}_0^*{\gz}_\alpha)^{-1}\\
&=(I_m-\lambda A)(I_m-\lambda QA)^{-1}{\gz}_\alpha({\gee}_0^*{\gz}_\alpha)^{-1}.
\end{split}
\]


\end{proof}

We turn next to the computation of $\langle (I-{\cX}-{\cY})^{-1}F{\gee}_0,F{\gee}_0\rangle_{{\cH }_G}$ via the recipe
\begin{equation}
\label{eq:jun14c24}
(I-{\cX}-{\cY})^{-1}=(I-{\cX})^{-1}\{I-{\cY}(I-{\cX})^{-1}\}^{-1}.
\end{equation}
The computation is carried out in steps under the assumption that $(H1)$--$(H3)$ are in force.
\bigskip

\noindent
{\bf 1.}\ {\it Verification of the formula}
\begin{equation}
\label{eq:mar8a24} 
(I-{\cX})^{-1}Fu=(I-b_\alpha(\omega)T_\alpha \Pi_\alpha)^{-1}Fu=F(I_m-XN_\alpha G)^{-1}u\quad for \ u\in\CC^m,
\end{equation}
{\it where }   $b_\alpha(\omega)=(\omega-\alpha)/(\omega-\ol{\alpha})=\rho_{\ol{\alpha}}(\omega)/\rho_\alpha(\omega)$,

\begin{equation}
\label{eq:jun14b24}
X=b_\alpha(\omega)(I_m-\ol{\alpha}A)(I_m-\alpha A)^{-1}
\quad \textrm{and}\quad (I_m-X)=\frac{\alpha-\ol{\alpha}}{\omega-\ol{\alpha}}(I_m-\omega A)(I_m-\alpha A)^{-1}.
\end{equation}
\bigskip

Let 
\[
(I-b_\alpha(\omega)T_\alpha \Pi_\alpha)^{-1}Fu=Fv.
\]
Then 
\[
\begin{split}
Fu&=(I-b_\alpha(\omega)T_\alpha \Pi_\alpha)Fv=F\{I_m-b_\alpha(\omega)(I_m-\ol{\alpha}A)(I_m-\alpha A)^{-1}N_\alpha G\}v\\ &=F(I_m-XN_\alpha G)v.
\end{split}
\]
Therefore, $v=(I_m-XN_\alpha G)^{-1}u$ and hence \eqref{eq:mar8a24} holds. 
\bigskip

\noindent
{\bf 2.}\ {\it Verification of the formula}
\begin{equation}
\label{eq:mar8b24}
\Pi_\alpha^\perp (I-b_\alpha(\omega)T_\alpha \Pi_\alpha)^{-1}Fu=F{\gz}_\alpha\{F(\omega){\gz}_\alpha\}^{-1}F(\omega)u.
\end{equation}
\bigskip

In view of \eqref{eq:mar8a24}, 
\[
\begin{split}
\Pi_\alpha^\perp(I-b_\alpha(\omega)T_\alpha \Pi_\alpha)^{-1}Fu&=F{\gz}_\alpha{\gz}_\alpha^*G(I_m-XN_\alpha G)^{-1}u\\
&=F{\gz}_\alpha{\gz}_\alpha^*G\{I_m+(I_m-X)^{-1}X{\gz}_\alpha{\gz}_\alpha^*G\}^{-1}(I_m-X)^{-1}u\\
&=F{\gz}_\alpha\{I_p+{\gz}_\alpha^*G(I_m-X)^{-1}X{\gz}_\alpha\}^{-1}{\gz}_\alpha^*G(I_m-X)^{-1}u\\
&=F{\gz}_\alpha\{{\gz}_\alpha^*G[I_m+(I_m-X)^{-1}X]{\gz}_\alpha\}^{-1}{\gz}_\alpha^*G(I_m-X)^{-1}u\\
&=F{\gz}_\alpha\{{\gz}_\alpha^*G(I_m-X)^{-1}{\gz}_\alpha\}^{-1}{\gz}_\alpha^*G(I_m-X)^{-1}u\\
\end{split}
\]
Formula \eqref{eq:mar8b24} drops out upon substituting 
\[
\begin{split}
{\gz}_\alpha^*G(I_m-X)^{-1}&=K_\alpha(\alpha)^{-1/2}F(\alpha)\Gamma G\{\rho_\alpha(\omega)\rho_\alpha(\alpha)^{-1}(I_m-\alpha A)(I_m-\omega A)^{-1}\}\\
&=\rho_\alpha(\omega)\rho_\alpha(\alpha)^{-1}K_\alpha(\alpha)^{-1/2}F(\omega)
\end{split}
\]
in the last expression.
\bigskip

\noindent
{\bf 3.}\ {\it Verification of the formula}
\begin{equation}
\label{eq:mar8c24}
\{I-{\cY}(I-{\cX})^{-1}\}^{-1}Fu=F\left\{I_m-\frac{\rho_{\ol{\alpha}}(\omega)}{\rho_{\alpha}(\alpha)^{1/2}}{\gz}_{\ol{\alpha}}Y(\omega)\right\}u\quad for\ u\in\CC^m,
\end{equation}
where
\begin{equation}
\label{eq:jun16a24}
Y(\omega)=S(\omega)(E_-(\omega)S(\omega)+E_+(\omega))^{-1}F(\omega)
\end{equation}
\bigskip

With the aid of the formula \eqref{eq:mar8b24}, it is readily checked that
\[
\{I_m-{\cY}(I_m-{\cX})^{-1}\}Fv=F\{I_m+b_\alpha(\omega){\gz}_{\ol{\alpha}}S(\omega)(F(\omega){\gz}_\alpha)^{-1}F(\omega)\}v
\]
and hence that 
\[
\begin{split}
\{I-{\cY}(I-{\cX})^{-1}\}^{-1}Fu&=F\{I_m+b_\alpha(\omega){\gz}_{\ol{\alpha}}S(\omega)(F(\omega){\gz}_\alpha)^{-1}F(\omega)\}^{-1}u\\
&=F(I_m-BC)^{-1}u,
\end{split}
\]
with $B=-b_\alpha(\omega){\gz}_{\ol{\alpha}}S(\omega)$ and $C=(F(\omega){\gz}_\alpha)^{-1}F(\omega)$. Thus, 
\[
\begin{split}
I_p-CB&=I_p+b_\alpha(\omega)(F(\omega){\gz}_\alpha)^{-1}F(\omega){\gz}_{\ol{\alpha}}S(\omega)\\
&=I_p+(\rho_\alpha(\omega)F(\omega){\gz}_\alpha)^{-1}\rho_{\ol{\alpha}}(\omega)F(\omega){\gz}_{\ol{\alpha}}S(\omega)\\
&=I_p+E_+(\omega)^{-1}E_-(\omega)S(\omega),
\end{split}
\]
is invertible for every point $\omega\in\CC_+$.  Consequently, $I_m-BC$ is invertible, and
\[
\begin{split}
(I_m-BC)^{-1}&=I_m+B(I_p-CB)^{-1}C\\
&=I_m-b_\alpha(\omega){\gz}_{\ol{\alpha}}S(\omega)(I_p+E_+(\omega)^{-1}E_-(\omega)S(\omega))^{-1}(F(\omega){\gz}_\alpha)^{-1}F(\omega)\\
&=I_m-\frac{\rho_{\ol{\alpha}}(\omega)}{\rho_\alpha(\alpha)^{1/2}}{\gz}_{\ol{\alpha}}S(\omega)(I_p+E_+(\omega)^{-1}E_-(\omega)S(\omega))^{-1}E_+(\omega)^{-1}F(\omega)\\
&=I_m-\frac{\rho_{\ol{\alpha}}(\omega)}{\rho_\alpha(\alpha)^{1/2}}{\gz}_{\ol{\alpha}}Y(\omega),
\end{split}
\]
which justifies \eqref{eq:mar8c24}.
\bigskip

\noindent
{\bf 4.}\ {\it Observe that in view of \eqref{eq:jun14c24}, \eqref{eq:mar8a24} and \eqref{eq:mar8c24}, 
$\langle (I-{\cX}-{\cY})^{-1}Fu, Fv\rangle_{{\cH}_G}=\textcircled{1}-\textcircled{2}$, where
\[
\begin{split}
\textcircled{1}&=\langle (I-{\cX})^{-1}Fu, Fv\rangle_{{\cH}_G}\\ &=\langle F(I_m-XN_\alpha G)^{-1}u, Fv\rangle_{{\cH}_G}
=v^*G(I_m-XN_\alpha G)^{-1}u
\end{split}
\]
and
\[
\begin{split}
\textcircled{2}&=\langle (I-{\cX})^{-1}F\frac{\rho_{\ol{\alpha}}(\omega)}{\rho_{\alpha}(\alpha)^{1/2}}{\gz}_{\ol{\alpha}}Y(\omega)u, Fv\rangle_{{\cH}_G}\\ 
&=\langle F(I_m-XN_\alpha G)^{-1}\frac{\rho_{\ol{\alpha}}(\omega)}{\rho_{\alpha}(\alpha)^{1/2}}{\gz}_{\ol{\alpha}}Y(\omega)u,Fv\rangle_{{\cH}_G}\\
&=v^*G(I_m-XN_\alpha G)^{-1} \frac{\rho_{\ol{\alpha}}(\omega)}{\rho_{\alpha}(\alpha)^{1/2}}{\gz}_{\ol{\alpha}}Y(\omega)u.
\end{split}
\]}

\bigskip

\noindent
{\bf 5.}\ {\it Verification of the formula}
\begin{equation}
\label{eq:mar9a24}
{\gee}_0^*G(I_m-XN_\alpha G)^{-1}{\gz}_{\ol{\alpha}}=-\frac{\rho_\alpha(\omega)}{\rho_\alpha(\alpha)^{1/2}}
(E _-^\#(\omega))^{-1}.
\end{equation}

Since,
\[
E_-^\#(\omega)=-\frac{\rho_\alpha(\omega)}{\rho_\alpha(\alpha)^{1/2}}{\gz}_{\ol{\alpha}}^*F^\#(\omega),
\]
it suffices to show that 
\[
{\gee}_0^*G(I_m-XN_\alpha G)^{-1}{\gz}_{\ol{\alpha}}{\gz}_{\ol{\alpha}}^*F^\#(\omega)=I_p.
\]
But,
\[
\begin{split}
(I_m-XN_\alpha G)^{-1}{\gz}_{\ol{\alpha}}{\gz}_{\ol{\alpha}}^*F^\#(\omega)&=(I_m-XN_\alpha G)^{-1}(\Gamma- N_{\ol{\alpha}})F^{\#}(\omega)\\
&=\Gamma (\Gamma-XN_\alpha )^{-1}(\Gamma- N_{\ol{\alpha}})F^{\#}(\omega)\\
&=\Gamma (\Gamma-XN_\alpha )^{-1}(\Gamma- XN_\alpha)F^{\#}(\omega)\\
&=\Gamma F^{\#}(\omega)
%
\end{split}
\]
since $XN_\alpha F^\#(\omega)=N_{\ol{\alpha}}F^\#(\omega)$, by (4) of Theorem 5.1. Therefore, 
\[
{\gee}_0^*G(I_m-XN_\alpha G)^{-1}{\gz}_{\ol{\alpha}}{\gz}_{\ol{\alpha}}^*F^\#(\omega)={\gee}_0^*G\Gamma F^\#(\omega)
=I_p,
\]
as claimed.
\bigskip

\noindent
{\bf 6.}\ {\it Verification of the formula}
\begin{equation}
\label{eq:may29b24}
\Phi_{\gE}(\omega)=E_+^\circ(\omega)E_+(\omega)^{-1}=\frac{2}{\rho_\alpha(\omega)}{\gee}_0^*G(I_m-\ol{\alpha}A)(I_m-\omega QA)^{-1}{\gz}_\alpha({\gee}_0^*{\gz}_\alpha)^{-1}
\end{equation}
with $Q=I_m-{\gz}_\alpha({\gee}_0^*{\gz}_\alpha)^{-1}{\gee}_0^*=Q^2$.
\bigskip

This is an easy consequence of \eqref{eq:nov19b23} and \eqref{eq:may20b24}.
\bigskip

\noindent
{\bf 7.}\ {\it Verify the formula} 
\[
{\gee}_0^*G(I_m-XN_\alpha G)^{-1}{\gee}_0=\frac{\rho_\alpha(\omega)}{\rho_\alpha(\alpha)}\left\{ h_0-\frac{\rho_{\ol{\alpha}}(\omega)}{2}\Phi_{\gE}(\omega)\right\}.
\]
\bigskip

By definition,
\[
\begin{split}
(I_m-XN_\alpha G)^{-1}&=(I_m-X+X{\gz}_\alpha{\gz}_\alpha^*G)^{-1}\\
&=(I_m-X)^{-1}\{I_m+X{\gz}_\alpha{\gz}_\alpha^*G(I_m-X)^{-1}\}^{-1}\\
&=(I_m-X)^{-1}\{I_m-CB\}^{-1}
\end{split}
\]
with $C=-X{\gz}_\alpha$ and   $B={\gz}_\alpha^*G(I_m-X)^{-1}$. Thus, as 
\[
\begin{split}
I_p-BC&=I_p+{\gz}_\alpha^*G(I_m-X)^{-1}X{\gz}_\alpha={\gz}_\alpha^*G(I_m-X)^{-1}{\gz}_\alpha\\
&=\frac{\rho_\alpha(\omega)}{\rho_\alpha(\alpha)}K_\alpha(\alpha)^{-1/2}F(\omega){\gz}_\alpha=\frac{1}{\rho_\alpha(\alpha)^{1/2}}K_\alpha(\alpha)^{-1/2}E_+(\omega)
\end{split}
\]
is invertible and $B=\rho_\alpha(\omega)\rho_\alpha(\alpha)^{-1}K_\alpha(\alpha)^{-1/2}F(\omega)$,
\[
\begin{split}
\{I_m-CB\}^{-1}&=I_m+C(I_p-BC)^{-1}B\\
&=I_m-X{\gz}_\alpha\{\rho_\alpha(\omega)\rho_\alpha(\alpha)^{-1}K_\alpha(\alpha)^{-1/2}F(\omega){\gz}_\alpha\}^{-1}
\rho_\alpha(\omega)\rho_\alpha(\alpha)^{-1}K_\alpha(\alpha)^{-1/2}F(\omega)\\
&=I_m-X{\gz}_\alpha(F(\omega){\gz}_\alpha)^{-1}F(\omega).
\end{split}
\]
Therefore,
\[
{\gee}_0^*G(I_m-XN_\alpha G)^{-1}{\gee}_0={\gee}_0^*G(I_m-X)^{-1}{\gee}_0-{\gee}_0^*G(I_m-X)^{-1}X{\gz}_\alpha (F(\omega){\gz}_\alpha)^{-1},
\]
which coincides with the asserted formula, since 
\[
(I_m-X)^{-1}{\gee}_0=\frac{\rho_\alpha(\omega)}{\rho_\alpha(\alpha)}{\gee}_0,\quad 
(I_m-X)^{-1}X=\frac{\rho_{\ol{\alpha}}(\omega)}{\rho_\alpha(\alpha)}(I_m-\omega A)^{-1}(I_m-\ol{\alpha}A)\quad\textrm{and} \quad \Phi_{\gE}=E_+^\circ E_+^{-1}.
\]
\bigskip

\noindent
{\bf 8.}\ {\it Verify the formula}
\begin{equation}
\label{eq:mar15b24}
\begin{split}
\{(I-{\cX}-{\cY})^{-1}&F{\gee}_0,F{\gee}_0\rangle_{{\cH}_G}=
\frac{\rho_\alpha(\omega)}{\rho_\alpha(\alpha)}h_0\\
&-\frac{\rho_\alpha(\omega)\rho_{\ol{\alpha}}(\omega)}{2\rho_\alpha(\alpha)}\{\Phi_{\gE}(\omega)-
2(E_-^\#(\omega))^{-1}S(\omega)(E_-(\omega)S(\omega)+E_+(\omega))^{-1}\},
\end{split}
\end{equation}
by combining the formulas in Steps 4, 5 and 7.
\bigskip

\begin{thm}
\label{thm:jun5a24}
If $(H1)$--$(H3)$ are in force,
 then the formula 
\begin{equation}
\label{eq:jun5a24}
\left\{\frac{1}{\pi i}\int_{-\infty}^{\infty} \frac{1}{\mu-\omega}d\sigma(\mu):\ \sigma\in\Sigma_{\rm H} \right\}
=\{(E_-^\circ S+E_+^\circ)(E_-S+E_+)^{-1}:\ S\in{\cS}_{\gE}^{\ptp}(\CC_+)\}
\end{equation}
establishes a $1:1$ correspondence between the solutions $\sigma$ of the truncated matrix Hamburger moment problem and the set of mvf's $S\in{\cS}_{\gE}^{\ptp}(\CC_+)$.
\end{thm}

\begin{proof} Formula \eqref{eq:jun5a24} is obtained by combining formulas \eqref{eq:mar15b24} and \eqref{eq:jun4b24}, since
\[
\Phi_{\gE}-
2(E_-^\#(\omega))^{-1}S(\omega)(E_-(\omega)S(\omega)+E_+(\omega))^{-1}=(E_-^\circ S+E_+^\circ)(E_-S+E_+)^{-1}.
\]
The fact that the correspondence is $1:1$ follows from Theorem \ref{thm:jun4a24}.
\end{proof}
\begin{rem}
\label{rem:jn14a24}
The formulas in \eqref{eq:dec1b23} for $E_\pm$  
involve many block columns of the matrix $\Gamma=G^{-1}$, since  $\alpha\in\CC_+$ .  The space ${\cH}_G$ may be identified as a de Branges space ${\cB}(\gE)$ based on   a de Branges matrix $\gE$ that only involves the last two columns of $\Gamma$ as in \cite{d88} by choosing 
\[
E_+(\lambda)=\sqrt{\pi}F(\lambda) ((\lambda+i){\gu}_n-{\gw}_n)\quad\textrm{and}\quad 
E_-(\lambda)=\sqrt{\pi}F(\lambda) ((\lambda-i){\gu}_n-{\gw}_n),
\]
where
\begin{equation*}
{\gu}_n=\Gamma {\gee}_n (\gamma_{jj}^{(n)})^{-1/2}, \quad 
 \mathfrak{v} _n=\gee_{n-1}-\gee_n \{\gamma_{nn}^{(n)}\}^{-1}\gamma_{n,n-1}^{(n)}\quad \textrm{and}\quad {\gw}_n=
\Gamma{\gv}_n(\gamma_{nn}^{(n)})^{-1/2}. 
\end{equation*}
This pair is obtained by setting $\alpha=re^{i\theta}$ in \eqref{eq:dec1b23}  and letting $r\uparrow\infty$. The formulas  

	\begin{align}\label{eq: Han eq 1}
N_\bullet=\Gamma-{\gu}_n{\gu}_n^*,\quad		F(\lambda)\lambda N_\bullet=F(\lambda)A^\ast N_\bullet \quad\textrm{for every point $\lambda\in\CC$}
	\end{align}
and
	\begin{align}\label{eq:jun24c24}
		A^\ast N_\bullet-N_\bullet A+{\gw}_n {\gu}_n^\ast -{\gu}_n {\gw}_n^\ast =0
	\end{align}
are obtained as a byproduct of this analysis. (To verify the latter, 
	let $X$ denote the left hand side of \eqref{eq:jun24c24}. Then it is readily checked that 
	\[
		\gee_i^\ast GXG \gee_j = \gee_i^\ast G A^\ast \gee_j -\gee_i^\ast AG \gee_j 
		= \gee_i^\ast G\gee_{j+1}-\gee_{i+1}^\ast G\gee_j 
		=0
	\]
	for $i,j= 0, \cdots, n-1$, since $G$ is a block Hankel matrix.  Thus as
	\begin{align*}
	\label{eq:jun24c24}
		\gee_n^\ast X=\gee_{n-1}^\ast (\Gamma-{\gu}_n{\gu}_n^\ast )-\gee_{n-1}^\ast {\gu}_n {\gw}_n^\ast =0 
	\end{align*}
	and the matrix $\begin{bmatrix}
		G{\gee}_0 & \cdots & G \gee_{n-1} & \gee_n
	\end{bmatrix}=   G  \begin{bmatrix}
		{\gee}_0 & \cdots &  \gee_{n-1} & \Gamma \gee_n
	\end{bmatrix}$
	is invertible, $X=0$.)
	
Since ${\gee}_n^*N_\bullet=0$ and $AA^*+{\gee}_n{\gee}_n^*=I_m$, 	\eqref{eq:jun24c24} implies that 
\[
N_\bullet=AN_\bullet A+AB\quad\textrm{with}\ B={\gu}_n {\gw}_n^\ast -{\gw}_n {\gu}_n^\ast, 
\]
which, upon iteration, yields the following 
formula of the Gohberg-Heinig type 
\[
N_\bullet=\sum_{j=0}^nA^jBA^j\quad\textrm{i.e.,}\quad \Gamma={\gu}_n{\gu}_n^*+\sum_{j=0}^nA^j({\gu}_n {\gw}_n^\ast -{\gw}_n {\gu}_n^\ast)A^j.
\]

\end{rem}

\section{Epilogue}
\label{sec:epilogue}
The main theme of this paper is the 
 natural  connection between the two truncated moment problems under consideration and de Branges spaces. In both settings the set of solutions to the truncated matrix moment problem is parameterized by a linear fractional transformation 
 \[
 T_\Theta[S]=(\Theta_{11}S+\Theta_{12})(\Theta_{21}S+\Theta_{22})^{-1}\quad\textrm{based on the mvf}\quad
 \Theta(\lambda)=\frac{1}{\sqrt{2}}\begin{bmatrix}E_-^\circ&E_+^\circ\\ E_-&E_+
\end{bmatrix} 
 \]
 with entries $E_\pm$ and 
 \[
 E_+^\circ(\omega)=\Phi_{\gE}(\omega)E_+(\omega)\quad\textrm{and}
\quad E_-(\omega)=\Phi_{\gE}(\omega)E_-(\omega)-2(E_-^\#(\omega))^{-1}, 
 \]
 i.e., all the entries in $\Theta$ are based 
  on the entries $E_\pm$ in the de Branges matrix ${\gE}$ chosen to fit the setting (i.e., in accordance with assumptions $(T1)$--$(T3)$ or $(H1)$--$(H3)$) and the appropriate choice of $f^\#$ and $\Phi_{\gE}$ (which is based on $E_+$). 
  
 We  have developed the parameterization of the sets $\Sigma_T$ and $\Sigma_H$ via the theory of unitary extensions of isometric operators because of the central role played by the isometric operators $S_\alpha$ with $\alpha\in\DD$ and $T_\alpha$ with $\alpha\in\CC_+$. However, since there is a $1:1$ correspondence between the normalized solutions $\sigma\in\Sigma_T$ and the set of $\Phi$ in the Carath\'{e}odory class ${\cC}^{\ptp}(\DD)$ of the form
 \[
\{\Phi\in{\cC}^{\ptp}(\DD):\ 
  \Phi(\lambda)=h_0+2\lambda h_1+\cdots+2\lambda^n h_n+\cdots\}
 \]
 and analogously (in view of Theorem 4.7 in \cite{bro71}) between the set of normalized $\sigma\in\Sigma_H$ and 
  \[
\{\Phi\in{\cC}^{\ptp}(\CC_+):\  \Phi(\lambda)=h_0+\lambda^{-1} h_1+\cdots+\lambda^{-2n} h_{2n}+\cdots\quad\textrm{and $\Vert \nu\Phi(i\nu)\Vert$ is bounded on $(1,\infty)$}\}, 
 \]
 the truncated matrix moment problem can also be viewed as an interpolation problem (in which the objective is to describe the set of $\Phi$ in the  Carath\'{e}odory class 
 ${\cC}^{\ptp}(\Omega_+)$ which match $\Phi_{\gE}$ in a prescribed way). This too will lead to the mvf $\Theta$ described above.

The mvf $\Theta(\lambda)$ has very special structure:

 \begin{thm}
 \label{thm:jun10a24}
  If either $(T1)$--$(T3)$ or  $(H1)$--$(H3)$ are in force, and if 
  \begin{equation}
\Theta(\lambda)=\frac{1}{\sqrt{2}}\begin{bmatrix}E_-^\circ&E_+^\circ\\ E_-&E_+
\end{bmatrix} \quad and\quad \FF(\lambda)=\frac{1}{\sqrt{2}}\begin{bmatrix}\Phi_{\gE}F-2\Pi_+\Delta_{\gE}F\\ F(\lambda),\end{bmatrix}
\end{equation}	
	then
			\begin{align}\label{thm:dec22a23}
		J_p- 
		\Theta(\lambda)j_p 
		\Theta(\omega)^\ast=\rho_\omega(\lambda)\FF(\lambda)\Gamma \FF(\omega)^\ast\quad\text{for all points}\ \lambda,\omega\in\CC
	\end{align}
and hence
the $2p\times 2p$ matrix polynomial $\Theta(\omega)$ is 
 subject to the constraints 
\begin{align}\label{eq:october31f23}
	\Theta(\omega )j_p \Theta(\omega)^\ast \prec J_p \quad\text{for } \omega \in \Omega_+\quad\textrm{and} \quad 
\Theta(\omega )j_p \Theta(\omega)^\ast =J_p \quad\text{for } \omega \in \Omega_0,
\end{align}
where $\Omega_+$ stands for either $\CC_+$ or $\DD$, and  $\Omega_0$ denotes the boundary of $\Omega_+$, i.e.,  
 $\RR$ if $\Omega_+=\CC_+$ and $\TT$ if $\Omega_+=\DD$. 
 \end{thm}	

\begin{proof} The proof is by direct calculation. 
\end{proof}
 
 \begin{rem}
 \label{rem:jun20a24}
The mvf $\FF(\lambda)$ considered in Theorem  \ref{thm:jun10a24} is of the form 
 $\begin{bmatrix}C_1\\C_2\end{bmatrix}(I_m-\lambda A)^{-1}$. Moreover, 
 if $(T1)$--$(T3)$ are in force, then 
\begin{equation}
\label{eq:jun20a24}
C_1=-\frac{1}{\sqrt{2}}{\gee}_0^*\LL^*, \quad C_2=\frac{1}{\sqrt{2}}{\gee}_0^*\quad\textrm{and}\quad G-A^*GA=2\pi\begin{bmatrix}C_1^*&C_2^*\end{bmatrix}J_p\begin{bmatrix}C_1\\ C_2\end{bmatrix},
\end{equation}
whereas, 
 if $(H1)$--$(H3)$ are in force, then 
\begin{equation}
\label{eq:jun20b24}
C_1=-\frac{1}{\sqrt{2}\pi i}{\gee}_0^*G, \quad C_2=\frac{1}{\sqrt{2}}{\gee}_0^*\quad\textrm{and}\quad A^*G-GA=\begin{bmatrix}C_1^*&C_2^*\end{bmatrix}J_p\begin{bmatrix}C_1\\ C_2\end{bmatrix}.
\end{equation}
In both settings the space ${\cM}=\{\FF u:\,u\in\CC^p\}$ endowed with the inner product $\langle \FF u,\FF v\rangle_{{\cM}_G}=v^*G u$ is a RKHS of entire vvf's with RK $K_\omega^{\cM}(\lambda)=\FF(\lambda)\Gamma\FF(\omega)^*$, which  is invariant with respect to $R_\omega$ for every point $\omega\in\CC$.
The identity in \eqref{eq:jun20b24} is equivalent to the  fact that  
\begin{equation}
\label{eq:jun23a24}
\langle R_\alpha f,g\rangle_{{\cH}_G}-\langle f, R_\beta g\rangle_{{\cH}_G}-(\alpha-\ol{\beta})\langle R_\alpha f, R_\beta g\rangle_{{\cH}_G}=2\pi ig(\beta)^*J_pf(\alpha)
\end{equation}
for every choice of $\alpha,\beta\in\CC$ and $f,g\in{\cH}_G$, 
which is due to de Branges (see e.g., Section 5.3 in \cite{ard08}), whereas the identity  in \eqref{eq:jun20a24} is equivalent to the  fact that 
\begin{equation}
\label{eq:jun23b24}
\langle (I+\alpha R_\alpha)f, (I+\beta R_\beta)g\rangle_{{\cH}_G}
-\langle R_\alpha f, R_\beta g\rangle_{{\cH}_G}=g(\beta)^*J_pf(\alpha)
\end{equation}
for every choice of $\alpha,\beta\in\CC$ and $f,g\in{\cH}_G$, which is in Ball \cite{ba75}, though he credits it to Rovnyak. 
\end{rem}

In view of  \eqref{eq:october31f23}, the mvf $\Theta(\lambda)\gV$ belongs to the class of $J_p$-inner mvf's. This has important implications, a number of which are listed in the next theorem.

\begin{thm}
\label{thm:nov5b23}
	If $\Theta$ meets the constraint \eqref{eq:october31f23} and $S\in \mathcal{S}^{p\times p}(\Omega_+)$, then:
	\begin{enumerate}
	\item $\Theta(\omega )j_p \Theta^\#(\omega)=J_p$ and $\Theta(\omega)^\#J_p\Theta(\omega)= j_p$ 
for all points  $\omega\in\CC$ if $\Omega_+=\CC_+$ and for all points $\omega\in\CC\setminus\{0\}$ if $\Omega_+=\DD$.
		\item $\Theta_{21}(\omega) S(\omega)+\Theta_{22}(\omega)$ is invertible for every point $\omega \in \Omega_+$.  
		\item $\Theta_{22}(\omega)$ and $\Theta_{21}(\omega)$ are  invertible for every point $\omega\in\Omega_0$.
		\item $\Theta_{22}^{-1}\Theta_{21}$ belongs to the class ${\cS}_{\rm in}^{\ptp}(\Omega_+)$  of $\ptp$ inner 
		mvf 's with respect to $\Omega_+$. 
		\item The linear fractional  transformation 
\begin{align}\label{eq: expression for T Theta}
		T_\Theta [S]= (\Theta_{11}S+\Theta_{12}) (\Theta_{21}S+\Theta_{22})^{-1},
	\end{align}
	maps $S\in \mathcal{S}^{p\times p}(\Omega_+)$ into $\mathcal{C}^{p\times p}(\Omega_+)$.
\item If $\Phi \in \mathcal{C}^{p\times p}(\Omega _+)$, $X=\begin{bmatrix}
		-I_p & \Phi
	\end{bmatrix}$, then
		\begin{align}\label{eq: condition on X and Theta}
		X(\omega)\Theta(\omega)j_p \Theta(\omega)^\ast X(\omega)^\ast \succeq 0
	\end{align}
	for every point $\omega \in \Omega_+$ if and only if
	$\Phi=T_\Theta[S]$ for some mvf $S\in \mathcal{S}^{p\times p}(\Omega _+)$.
	\end{enumerate}

\end{thm}
\begin{proof}
See e.g.,  \cite{ard08} for $\Omega_+=\CC_+$; the proofs for $\Omega_+=\DD$ are much the same. 
\end{proof}

The supplementary identity
\begin{equation}
\label{eq:jun19a24}
T_\Theta[S]+T_\Theta[S]^*=(S^*\Theta_{21}^*+\Theta_{22}^*)^{-1}(I_p-S^*S)(\Theta_{21}S+\Theta_{22})^{-1}, 
\quad\textrm{a.e. on $\Omega_0$}
\end{equation}
which is obtained by straightforward calculation underlies maximum entropy inequalities:

{\it 
If $(H1)$--$(H3)$ are in force and $\Phi=T_\Theta[S]$, then for each point $\omega\in\CC_+$ and each mvf $S\in{\cS}^{\ptp}(\CC_+)$
\begin{equation}
\label{eq:nov22a23}
\frac{\omega-\ol{\omega}}{2\pi i}\int_{-\infty}^\infty \ln\,\det\,\{\Phi(\mu)+\Phi(\mu)^*\}\, \frac{d\mu}{\vert\mu-\omega\vert^2}\le 
\ln \det\{E_+(\omega)E_+(\omega)^*-E_-(\omega) E_-(\omega)^*\}
\end{equation}
with equality if and only if  $S(\lambda)\equiv-(E_+^{-1}E_-)(\omega)^*$. 

If  $(T1)$--$(T3)$ are in force and $\Phi=T_\Theta[S]$, then for each point $\omega\in\DD$ and each mvf $S\in{\cS}^{\ptp}(\DD)$
\begin{equation}
\label{eq:nov22aa23}
\frac{1-\vert \omega\vert^2}{2\pi}\int_0^{2\pi} \ln\,\det\,(\Phi(e^{it})+\Phi(e^{it})^*)\, \frac{dt}{\vert e^{it}-\omega\vert^2}d\mu\le 
\ln \det\{E_+(\omega)E_+(\omega)^*-E_-(\omega)E_-(\omega)^*\}
\end{equation}
with equality if and only if $S(\lambda)\equiv-(E_+^{-1}E_-)(\omega)^*$.}

These inequalities  are discussed at assorted levels of generality in  \cite{d89a}, 
Section 11 of \cite{d89b}, \cite{ga92}, \cite{mm05}, Chapter 11 of \cite{ard08} (which is based largely on \cite{ark81} and \cite{ark83})  and a number of the references cited therein. 

\vspace{1cm}

\noindent \textbf{Acknowledgment}
 The authors thank Professor Vladimir Derkach for sharing his expertise on the application of selfadjoint extensions of symmetric operators and unitary extensions of isometric operators to moment problems. K. Dhara thanks the Weizmann Institute of Science
	for a postdoctoral fellowship and
	the Department of Science and Technology, Govt. of India, for  financial support via an INSPIRE Faculty Fellowship (Ref no: DST/INSPIRE/04/2022/003288). The authors also thank the referee for useful remarks. 
 
\noindent \textbf{Competing interests } 
The authors declare that there is no conflict of interest.

\noindent \textbf{Data Availability}
Not applicable.

\begin{thebibliography}{99999}
	
	
	
	\bibitem[Ak65]{ak65}
N. I. Akhiezer, {\it The Classical Moment Problem}, Oliver and Boyd, Edinburgh and London, 1965.







\bibitem[An70]{an70}	
T. Ando, {\it Truncated moment problems for operators}, 
Acta Sci. Math. (Szeged),  {\bf 31} (1970), 319–334.




\bibitem[ArD08]{ard08}
D. Z. Arov and H. Dym, {\it $J$-contractive Matrix Valued Functions and
	Related Topics}, Cambridge University Press,
Cambridge, 2008.




\bibitem[ArD18]{ard18}
D. Z. Arov and H. Dym, {\it Multivariate Prediction, de Branges Spaces, and Related Extension and Inverse Problems}, Birkhauser, Basel, 2018.


\bibitem[Ark81]{ark81}
D. Z. Arov and M. G. Krein, {\it The problem of finding the minimum entropy in indeterminate problems of continuation}, (Russian) Funktsional Anal. i Prilozhen, {\bf 15} (1981), no. 2, 61--64.

\bibitem[Ark83]{ark83}
D. Z. Arov and M. G. Krein, {\it Calculation of entropy functionals and their minima in indeterminate continuation problems}, (Russian) Acta Sci. Math. (Szeged) {\bf 45} (1983), 33-50.

\bibitem[Ba75]{ba75} 
J. A. Ball, {\it Models for noncontractions}, J. Math. Anal. Appl. {\bf 52} (1975), 235--254.

\bibitem[BD96]{bd96}
C. Berg and A. J. Duran, {\it Orthogonal polynomials, $L_2$ spaces and entire functions}, Math. Scand., {\bf 79} (1996), 209--223.


\bibitem[Br63]{br63}
L. de Branges,
{\it Some Hilbert spaces of analytic functions I}, Trans. Amer. Math. Soc.,
{\bf 106} (1963), 445--668.
\bibitem[Br68a]{br68a}
L. de Branges, {\it Hilbert Spaces of Entire Functions}, Prentice Hall, London, 1968.


\bibitem[Br68b]{br68b}
L.  de Branges,  {\it The expansion theorem for Hilbert spaces of entire functions}, in:   Entire Functions and Related Parts of Analysis, Proc. Sympos. Pure Math., American Mathematical Society, Providence, Rhode Island, 1968,  pp. 79--148.

\bibitem[BrR66]{brr66}
L. de Branges and J. Rovnyak,
{\it Canonical models in quantum scattering theory},
in: Perturbation Theory and its Application in Quantum Mechanics,
Wiley, New York, 1966, pp. 359--391.

\bibitem[Bro71]{bro71}
M. S. Brodskii, {\it Triangular and Jordan representation of Linear Operators}, Translations of Mathematical Monographs, vol. 32, American Mathematical Society, 1971.

\bibitem[Ch67]{ch67}
M. E. Chumakin, {\it Generalized resolvents of isometric operators}, 	Sibirsk. Mat. Ž. 8 (1967), 876–892.

\bibitem[D88]{d88}
H. Dym, {\it Hermitian block Toeplitz matrices, orthogonal polynomials, reproducing kernel Pontryagin spaces, interpolation and extension},
Oper. Theory Adv. Appl., 34
Birkhäuser Verlag, Basel, 1988.


\bibitem[D89a]{d89a}
H.	Dym, {\it 
	On Hermitian block Hankel matrices, matrix polynomials, the Hamburger moment problem, interpolation and maximum entropy}, 
Integral Equations Operator Theory {\bf 12}(1989), no.6, 757–812.


\bibitem[D89b]{d89b}
H. Dym, {\it  $J$-contractive matrix functions, reproducing kernel Hilbert spaces and interpolation}, CBMS Regional Conference Series in Mathematics, {\bf 71}. Published for the Conference Board of the Mathematical Sciences, Washington, DC; by the American Mathematical Society, Providence, Rhode Island, 1989.



\bibitem[D23a]{d23a}
H. Dym, {\it Two classes of vector valued de Branges spaces}, J. Functional Analysis, {\bf 284} (2023), 109758


\bibitem[D23b]{d23b}
H. Dym, {\it Linear Algebra in Action}, third edition, American Mathematical Society, Providence, Rhode Island, 2023.


\bibitem[DS17]{ds17}
H. Dym and S. Sarkar,  {\it Multiplication operators with deficiency indices $(p,p)$ and sampling formulas in reproducing kernel Hilbert spaces of entire vector valued functions}, J. Functional Analysis, {\bf 273} (2017), 3671--3718.


\bibitem[Ga92]{ga92}
J. P. Gabardo, {\it A maximum entropy approach to the classical moment problem}, J.  Functional Analysis, {\bf 106}
(1992), 80--94.

\bibitem[GH74]{gh74}
I. C. Gohberg and G. Heinig, {\it Inversion of finite Toeplitz matrices consisting of elements of a noncommutative algebra}, (Russian)
Rev. Roumaine Math. Pures Appl., 19 (1974), 623–663.

\bibitem[MM03]{mm03}
J. G. Marcano and M. D. Morán, {\it 
	The Arov-Grossman model and the Burg multivariate entropy}, 
J. Fourier Anal. Appl. {\bf 9} (2003), no.6, 623–647.

\bibitem[MM05]{mm05}
J. G. Marcano and M. D. Morán,
{\it The Arov-Grossman model and Burg's entropy}, Recent advances in applied probability, 329-349.
Springer-Verlag, New York, 2005.



\bibitem[Za13]{za13}
S. M. Zagorodnyuk, 
Nevanlinna formula for the truncated matrix trigonometric moment problem,
Ukrainian Math. J., {\bf 64} (2013), no. 8, 1199–1214.


\bibitem[Za12]{za12}	
S. M. Zagorodnyuk, The Nevanlinna-type formula for the matrix Hamburger moment problem, 
Methods Funct. Anal. Topology, {\bf 18} (2012), no. 4, 387–400.

\bibitem[Za11]{za11}	
S. M. Zagorodnyuk, Truncated matrix trigonometric problem of moments: operator approach, 
Ukrainian Math. J., {\bf 63} (2011), no. 6, 914–926.





\end{thebibliography}
\end{document}